\documentclass{amsart}
\usepackage{times}
\usepackage{amsfonts,amssymb,amscd,txfonts}
\usepackage{graphicx}

\newtheorem{thm}{Theorem}
\newtheorem{defi}{Definition}
\newtheorem{lem}{Lemma}
\newtheorem{cor}{Corollary}
\newtheorem{pro}{Proposition}
\newtheorem{rem}{Remark}

\newtheorem{obs}{Observation}
\newcommand{\pf}{{\bf\sl Proof.  }}
\renewcommand{\hbar}{\centerline{\rule{8cm}{0.5mm}}}
\renewcommand{\qed}{\vrule width0pt\hfill \raisebox{-.3ex}
   {\vrule height8pt width8pt depth0pt} \hspace*{-7pt}}
\renewcommand{\Re}{{\mathbb R}}
\newcommand{\Z}{{\mathbb Z}}
\newcommand{\C}{{\mathbb C}}

\newcommand{\ctc}{{\rm CTC}}
\newcommand{\tc}{{\rm TC}}

\newcommand{\nc}{{\rm NTC}}

\begin{document}

\large

\title[Total Curvature of Graphs after Milnor and Euler]
{\textbf {Total Curvature of Graphs after Milnor and Euler}}
\author{\bf Robert Gulliver and Sumio Yamada}
\thanks{Supported in part by JSPS Grant-in-aid for Scientific 
Research No.17740030}
\thanks{Thanks to the Korea Institute for Advanced
Study for invitations.}
\date{December 31, 2010}

\begin{abstract}
We define a new notion of total curvature, called {\em net total
curvature,} for finite graphs
embedded in $\Re^n$, and investigate its properties.  Two guiding
principles are given by Milnor's way of measuring the local 
crookedness
of a Jordan curve via a Crofton-type formula, and by considering
the double cover of a given graph as an Eulerian circuit.  The
strength of combining these ideas in defining the curvature
functional is (1) it allows us to interpret the
singular/non-eulidean behavior at the vertices of the graph as a
superposition of vertices of a $1$-dimensional manifold, and thus
(2) one can compute the total curvature for a wide range of graphs
by contrasting local and global properties of the graph utilizing
the integral geometric representation of the curvature.  A
collection of results on upper/lower bounds of the total curvature
on isotopy/homeomorphism classes of embeddings is presented, which
in turn demonstrates the effectiveness of net total curvature as a
new functional 

measuring complexity of spatial
graphs in differential-geometric terms.  
\end{abstract}

\maketitle

\section{INTRODUCTION: CURVATURE OF A GRAPH}\label{intro}

	The celebrated F\'ary-Milnor theorem states that a curve in
$\Re^n$ of total curvature at most $4\pi$ is unknotted.  

As a key step in his 1950 proof, John Milnor showed that for a
smooth
Jordan curve $\Gamma$ in $\Re^3$, the total curvature equals
half the integral over $e \in S^2$ of the number $\mu(e)$ of
local maxima of the linear ``height" function $\langle e,\cdot \rangle$
along $\Gamma$ \cite{M}.  This equality can be regarded as a
Crofton-type representation formula of total curvature where the
order of integrations over the curve and the unit tangent sphere
(the space of directions) are reversed.  The F\'ary-Milnor theorem
follows, since total curvature less than $4\pi$ implies there is a
unit vector $e_0 \in S^2$ so that $\langle e_0,\cdot \rangle$ has
a unique local maximum, and therefore that this linear function
is increasing on an interval of $\Gamma$ and decreasing
on the complement.  Without changing the pointwise value of this
``height" function, $\Gamma$ can be topologically untwisted to a
standard embedding of $S^1$ into $\Re^3$.  The Fenchel theorem,
that any curve in $\Re^3$ has total curvature at least $2\pi$,
also follows from Milnor's key step, since for all $e\in S^2$, the
linear function $\langle e,\cdot \rangle$ assumes its maximum
somewhere along $\Gamma$, implying $\mu(e) \geq 1$.  Milnor's
proof is independent of the proof of Istvan F\'ary, published
earlier, which takes a different approach \cite{Fa}.

       We would like to extend the methods of Milnor's seminal
paper, replacing the
simple closed curve by a finite {\em graph} $\Gamma$ in $\Re^3$.
$\Gamma$ consists of a finite number of points, the
{\em vertices}, and a finite number of simple arcs, the
{\em edges}, each of which has as its endpoints one or two of 
the vertices.  We
shall assume $\Gamma$ is connected.  The {\em degree} of a vertex
$q$ is the number $d(q)$ of edges which have $q$ as an endpoint.
(Another word for degree is ``valence".) We remark that it is
technically not needed that the dimension $n$ of the ambient space
equals three.  All the arguments can be generalized to higher
dimensions, although in higher dimensions 
$(n\geq 4)$ there are no nontrivial knots. Moreover, any two 	
homeomorphic 		
graphs are isotopic.   

The key idea in generalizing total curvature for knots to total
curvature for graphs is to consider the Euler circuits of the
given graph, namely, parameterizations by $S^1$, of
the {\it double} cover of the graph. We note that given 
a graph of even degree, there can be several Euler circuits, or
ways to ``trace it without
lifting the pen." A topological vertex of a graph of degree $d$
is a singularity, in that the graph is not
locally Euclidean.  However by considering an Euler circuit of
the double of the graph, the vertex 
becomes locally the intersection point of $d$ paths. We
will show (Corollary~\ref{cor2}) that at the vertex, each path
through it has a (signed) measure-valued curvature, and
the  absolute value of the sum of those measures is well-defined,
independent of the choice of the Euler circuit of the double
cover. We define (Definition~\ref{defnet}) the {\em net total
curvature} (NTC) of a piecewise $C^2$ graph to be the sum of the
total curvature of the smooth arcs and the contributions from the
vertices as described.

This notion of net total curvature is substantially different from
the total curvature, denoted TC, as defined by Taniyama \cite{T}.
(Taniyama writes $\tau$ for TC.) See section \ref{deftc} below.

This is consistent with known results for the vertices of degree
$d =2$; with vertices of degree three or more, this definition
helps facilitate a new Crofton-type representation formula
(Theorem~\ref{muthm}) for total curvature of graphs, where the
total curvature is represented as an integral over the unit
sphere. Recall that the vertex is now seen as $d$ distinct points
on an Euler circuit. The way we pick up the contribution of the
total curvature at the vertices identifies the $d$ distinct
points, and thus the $2d$ unit tangent spheres on a circuit.  As
Crofton's formula in effect reverses the order of integrations ---
one over the circuit, the other over the space of tangent
directions --- the sum of the $d$ exterior angles at the vertex is
incorporated in the integral over the unit sphere.  On the other
hand the integrand of the integral over the unit sphere counts the
number of net local maxima of the height function along an axis,
where net local maximum means the number of local maxima minus the
number of local minima at these $d$ points of the Euler circuit.
This establishes a correspondence between the differential
geometric quantity (net total curvature) and the differential
topological quantity 
(average number of maxima) of the graph, as stated in 
Theorem~\ref{muthm} below. 

In section~\ref{deftc}, we compare several definitions for total
curvature of graphs which have appeared in the recent literature. 
In section~\ref{comb}, we introduce the main tool (Lemma
~\ref{combin}) which in a sense reduces the computation of NTC
to counting intersections with planes.

Milnor's treatment~\cite{M} of total curvature also contained an
important topological extension.  Namely, in order to define total
curvature, the knot needs only to be {\em continuous}.  This makes
the total curvature a geometric quantity defined on any
homeomorphic image of $S^1$.  In this article, we first define net
total curvature (Definition~\ref{defnet}) on piecewise $C^2$
graphs, and then extend the definition to continuous graphs
(Definition~\ref{gendefnet}.)  In analogy to Milnor, we
approximate a given continuous graph by a sequence of polygonal
graphs.  In showing the monotonicity of the total curvature
(Proposition~\ref{monotmu}) under the refining process of
approximating graphs we use our representation formula
(Theorem~\ref{muthm}) applied to the polygonal graphs.  

Consequently the Crofton-type representation formula
is also extended (Theorem~\ref{muthm2}) to cover continuous
graphs.  Additionally, we are able to show that continuous 
graphs with finite total curvature 
(NTC or TC) 
are tame.  We say that a graph is {\it tame} when
it is isotopic to an embedded polyhedral graph.

In sections~\ref{three/four} through~\ref{FaryMilnor}, we characterize 
NTC with respect to the geometry and the topology of the graph.  
Proposition~\ref{subadditivity} shows the subadditivity of NTC under
the union of graphs which meet in a finite set.  In section
~\ref{deg3}, the concept of bridge number is extended from knots to
graphs, in terms of which the minimum of NTC can be explicitly
computed, provided the graph has at most one vertex of degree $>3$.
In section~\ref{lowbds}, Theorem~\ref{incrdecr} gives a lower bound 
for NTC in terms of the width of an isotopy class. The infimum of
NTC is computed for specific graph types: the two-vertex graphs
$\theta_m$, the ``ladder" $L_m$, the ``wheel" $W_m$, the complete
graph $K_m$ on $m$ vertices and the complete bipartite graph $K_{m,n}$.

Finally we prove a result
(Theorem~\ref{thetathm})
which gives a
Fenchel type lower bound $(\geq 3 \pi)$ for total curvature of a
theta graph (an image of the graph consisting of a circle with an
arc connecting a pair of antipodal points), and a F\'{a}ry-Milnor
type upper bound $(< 4 \pi)$ to imply the theta graph is isotopic
to the standard embedding. A similar result was given by Taniyama
\cite{T}, referring to TC. In contrast, for graphs of the type of
$K_m \ (m\geq 4)$, the infimum of NTC in the 
isotopy class of 
a polygon on $m$ vertices is also the infimum for a 
sequence of distinct isotopy classes. 

Many of the results in our earlier preprint \cite{GY2} have been
incorporated into the present paper.

We thank Yuya Koda for his comments regarding
Proposition~\ref{net3}, and Jaigyoung Choe and Rob Kusner for
their comments about Theorem~\ref{thetathm}, especially about the
sharp case $\nc(\Gamma) = 3\pi$ of the lower bound estimate.

\section{DEFINITIONS OF TOTAL CURVATURE}\label{deftc}

\vspace{2mm}

 The first difficulty,
in extending the results of Milnor's classic paper, is to
understand the contribution to total curvature at a vertex of
degree $d(q)\geq 3$. We first consider the well-known case:

\vspace{2mm}

\centerline{\bf{Definition of Total Curvature for Knots}}

\vspace{2mm}
      
For a smooth closed curve $\Gamma$, the total curvature is
$$
{\mathcal C}(\Gamma) = \int_\Gamma |\vec{k}| \, ds,
$$
where $s$ denotes arc length along $\Gamma$ and $\vec{k}$ is the
curvature vector.  If $x(s)\in \Re^3$ denotes the position of the
point measured at arc length $s$ along the curve, then
$\vec{k} = \frac {d^2x}{ds^2}$.  For a piecewise smooth curve,
that is, a graph with vertices $q_1, \dots, q_N$ having always
degree $d(q_i)=2$, the total curvature is readily generalized to

\begin{equation}\label{gencurv}
{\mathcal C}(\Gamma) =
\sum_{i=1}^N {\rm c}(q_i) + \int_{\Gamma_{\rm reg}} |\vec{k}| \,
ds,
\end{equation}
where the integral is taken over the separate $C^2$ edges of
$\Gamma$ without their endpoints;
and where ${\rm c}(q_i) \in [0,\pi]$ is the
exterior angle formed by the two
edges of $\Gamma$ which meet at $q_i$.  That is,
$\cos({\rm c}(q_i)) = \langle T_1, -T_2\rangle,$
where $T_1= \frac{dx}{ds}(q_i^+)$ and
$T_2= -\frac{dx}{ds}(q_i^-)$ are the unit tangent vectors at
$q_i$ pointing into the two edges which meet at $q_i$.
The exterior angle ${\rm c}(q_i)$ is the correct contribution to
total
curvature, since any sequence of smooth curves converging to
$\Gamma$ in $C^0$, with $C^1$ convergence on compact subsets of
each open edge, includes a small arc near $q_i$ along which the
tangent vector changes from near $\frac{dx}{ds}(q_i^-)$ to near
$\frac{dx}{ds}(q_i^+)$.  The greatest lower bound of the
contribution to total curvature of this disappearing
arc along the smooth approximating curves equals ${\rm c}(q_i)$.

Note that ${\mathcal C}(\Gamma)$ is well defined for an {\em
immersed} knot $\Gamma$.

\vspace{2mm}

\centerline{\bf{Definitions of Total Curvature for Graphs}}

\vspace{2mm}

When we turn our attention to a {\em graph} $\Gamma$, we find the
above definition for curves (degree $d(q)=2$) does not generalize
in any obvious way to higher degree (see \cite{G}).
The ambiguity of the general formula \eqref{gencurv} is resolved
if we specify the replacement for ${\rm c}(0)$ when $\Gamma$ is
the cone over a finite set $\{T_1, \dots, T_d\}$ in the unit
sphere $S^2$.  

The earliest
notion of total curvature of a graph appears in the
context of the first variation of length of a graph, which we call
{\bf variational total curvature}, and is  
called the {\em mean curvature} of the graph in \cite{AA}: we
shall write VTC.  The contribution to VTC at a vertex $q$ of
degree $2$, with unit tangent vectors $T_1$ and $T_2$, is 
${\rm vtc}(q)=|T_1 + T_2| = 2\sin(c(q)/2)$. 
At a non-straight vertex $q$ of degree $2$, ${\rm vtc}(q)$ is less
than the exterior angle ${\rm c}(q)$.  For a vertex of degree $d$, the
contribution is ${\rm vtc}(q)=|T_1 +\dots +T_d|$.

A rather natural definition of total curvature of graphs was 
given by Taniyama in \cite{T}.  We have called this {\bf maximal 
total curvature} $\tc(\Gamma)$ in \cite{G}.  The contribution to 
total curvature at a vertex $q$ of degree $d$ is
$${\rm tc}(q):=
\sum_{1\leq i<j\leq d}\arccos\langle T_i,-T_j\rangle.$$
In the case $d(q) = 2$, the sum above has only one term, the
exterior angle ${\rm c}(q)$ at $q$.  Since the length of the Gauss
image of a curve in $S^2$ is the total curvature of the curve,
${\rm tc}(q)$ may be interpreted as adding to the Gauss image in
$\Re P^2$ of the edges, a complete great-circle graph on
$T_1(q),\dots,T_d(q)$, for each vertex $q$ of degree $d$. Note
that the edge between two vertices does not measure the distance
in $\Re P^2$ but its supplement.

In our earlier paper \cite{GY1} on the density of an
area-minimizing two-dimensional rectifiable set $\Sigma$ spanning
$\Gamma$, we found that it was very useful to apply the
Gauss-Bonnet formula to the cone over $\Gamma$ with a point $p$ of
$\Sigma$ as vertex.  The relevant notion of total curvature in
that context is {\bf cone total curvature} $\ctc(\Gamma)$, defined
using ${\rm ctc}(q)$ as the replacement for ${\rm c}(q)$ in
equation \eqref{gencurv}:

\begin{equation}\label{defconetc}
{\rm ctc}(q) := \sup_{e \in S^2} \left\{
\sum_{i=1}^d\left(\frac{\pi}{2}-\arccos\langle T_i, e\rangle
\right) \right\}.
\end{equation}
Note that in the case $d(q) = 2$, the supremum above is assumed at
vectors $e$ lying in the smaller angle between the tangent vectors
$T_1$ and $T_2$ to $\Gamma$, so that ${\rm ctc}(q)$ is then the
exterior angle ${\rm c}(q)$ at $q$.  The main result of \cite{GY1}
is that $2\pi$ times the area density of $\Sigma$ at any of its
points is at most equal to $\ctc(\Gamma)$.  The same result had
been proven by Eckholm, White and Wienholtz for the case of a
simple closed curve \cite{EWW}.  Taking $\Sigma$ to be the
branched immersion of the disk given by Douglas \cite{D1} and
Rad\'o \cite{R}, it follows that if ${\mathcal C}(\Gamma) \leq
4\pi$, then $\Sigma$ is embedded, and therefore $\Gamma$ is
unknotted.  Thus \cite{EWW} provided an independent proof of the
F\'ary-Milnor theorem.  However, $\ctc(\Gamma)$ may be small for
graphs which are far from the simplest isotopy types of a graph
$\Gamma$. 

In this paper, we introduce the notion of {\bf net total
curvature} $\nc(\Gamma)$, which is the appropriate definition for
generalizing --- {\em to graphs} --- Milnor's approach to isotopy
and total curvature of {\em curves}.  For each unit tangent vector
$T_i$ at $q$, $1 \leq i \leq d=d(q)$, let $\chi_i:S^2 \rightarrow
\{-1, +1\}$ be equal to $-1$ on the hemisphere with center at
$T_i$, and $+1$ on the opposite hemisphere (modulo sets of zero
Lebesgue measure).  We then define

\begin{equation}\label{defnc}
{\rm ntc}(q):=
\frac{1}{4}\int_{S^2}\left[\sum_{i=1}^d\chi_i(e)\right]^+\,dA_{S^2}(e).
\end{equation}
We note that the function $\sum_{i=1}^d \chi_i(e)$ is odd, hence 
the quantity above can be written as  
\[
{\rm ntc}(q):=
\frac{1}{8}\int_{S^2}\left|\sum_{i=1}^d\chi_i(e)\right|\,dA_{S^2}(e).
\]
as well.  In the case $d(q)=2$, the integrand of \eqref{defnc} is
positive (and equals 2) only on the set of unit vectors $e$ which
have negative inner products with both $T_1$ and $T_2$, ignoring
$e$ in sets of measure zero.  This set is bounded by semi-great 
circles orthogonal to $T_1$ and to $T_2$, and has spherical area
equal to twice 
the exterior angle.  So in this case, 
${\rm ntc}(q)$ is the exterior angle. Thus, in the special case
where $\Gamma$ is a piecewise smooth curve, the following quantity
$\nc(\Gamma)$ coincides with total curvature, as well as with
$\tc(\Gamma)$ and $\ctc(\Gamma)$:

%
%
\begin{defi}\label{defnet}
We define the {\em net total curvature} of a piecewise $C^2$
graph $\Gamma$ with vertices $\{q_1, \dots, q_N\}$ as
%
%
\begin{equation}
{\nc}(\Gamma):=
\sum_{i=1}^N {\rm ntc}(q_i)+\int_{\Gamma_{\rm reg}} |\vec{k}| \,
ds.
\end{equation}
\end{defi}

For the sake of simplicity, elsewhere in this paper, we consider
the ambient space to be ${\Re}^3$.  However the definition of the
net total curvature can be generalized for a graph in ${\Re}^n$ by
defining the vertex contribution in terms of an average over
$S^{n-1}$:
\[
{\rm ntc}(q) := 
\pi \Big(\fint_{S^{n-1}}\left[\sum_{i=1}^d
\chi_i(e)\right]^+ \,dA_{S^{n-1}}(e) \Big),
\]  
which is consistent with the definition~(\ref{defnc}) of $\rm ntc$
when $n=3$. 

Recall that Milnor~\cite{M} defines the total curvature of a
continuous simple closed curve $C$ as the supremum of the total
curvature of all polygons inscribed in $C$. By analogy, we define
net total curvature of a {\it continuous} graph $\Gamma$ to be the
supremum of the net total curvature of all polygonal graphs $P$
suitably inscribed in $\Gamma$ as follows.  

%
%
\begin{defi}\label{gamapprox}
For a given continuous graph $\Gamma$, we
say a polygonal graph $P \subset \Re^3$ is 
{\em $\Gamma$-approximating}, provided that its topological
vertices (those of degree $\neq 2$) are exactly the topological
vertices of \, $\Gamma$, and having the same degrees; and that the
arcs of $P$ between two topological vertices correspond one-to-one 
to the edges of \, $\Gamma$ between those two vertices.
\end{defi}
Note that if $P$ is a $\Gamma$-approximating polygonal graph, then
$P$ is homeomorphic to $\Gamma$.   According to the statement of
Proposition \ref{monotmu}, whose proof will be given in the next
section, if $P$ and $\widetilde{P}$ are $\Gamma$-approximating
polygonal graphs, and $\widetilde{P}$ is a refinement of $P$, then
$\nc(\widetilde{P}) \geq \nc(P)$. Here $\widetilde{P}$ is said to
be a refinement of $P$ provided the set of vertices of $P$ is a
subset of the vertices of $\widetilde{P}$.  Assuming Proposition
\ref{monotmu} for the moment, we can generalize the definition of
the total curvature to non-smooth graphs.

%
\begin{defi}\label{gendefnet}
Define the {\em net total curvature} of a continuous graph 
$\Gamma$ by
\[
\nc(\Gamma) := \sup_{P} \nc(P)
\]
where the supremum is taken over all $\Gamma$-approximating
polygonal graphs $P$.  
\end{defi}

For a polygonal graph $P$, applying Definition \ref{defnet}, 
$$
{\nc}(P):=
\sum_{i=1}^N {\rm ntc}(q_i), 
$$
where $q_1, \dots, q_N$ are the vertices of $P$.

Definition \ref{gendefnet} is consistent with Definition
\ref{defnet} in the 
case of a piecewise $C^2$ graph $\Gamma$.  Namely, as Milnor 
showed, the total curvature ${\mathcal C}(\Gamma_0)$ of 
a smooth curve $\Gamma_0$ is 
the supremum of the total curvature of inscribed polygons 
(\cite{M}, p. 251), which gives the required supremum for each
edge.  At a vertex $q$ of the piecewise-$C^2$ graph $\Gamma$, 
as a sequence $P_k$ of  $\Gamma$-approximating 
polygons become arbitrarily fine, a vertex $q$ of $P_k$
(and of $\Gamma$) has unit tangent vectors converging in $S^2$
to the unit tangent vectors to $\Gamma$ at $q$.  
It follows that for 
$1\leq i\leq d(q)$, $\chi_i^{P_k} \to \chi_i^{\Gamma}$ 
in measure on $S^2$, and therefore 
${\rm ntc}_{P_k}(q) \to {\rm ntc}_\Gamma(q)$.

%
\section{CROFTON-TYPE REPRESENTATION FORMULA FOR TOTAL
CURVATURE}\label{comb}
We would like to explain how the net total curvature $\nc(\Gamma)$ 
of a graph is related to more familiar notions of
total curvature.  Recall that a graph $\Gamma$ has an Euler
circuit if and only if its vertices all have even degree, by a
theorem of Euler.  An Euler circuit is a closed, connected path
which traverses each edge of $\Gamma$ exactly once.  Of course, we
do not have the hypothesis of even degree.  We can attain that
hypothesis by passing to the {\em double} $\widetilde{\Gamma}$ of
$\Gamma$:  $\widetilde{\Gamma}$ is the graph with the same
vertices as $\Gamma$, but with two copies of each edge of
$\Gamma$.  Then at each vertex $q$, the degree as a vertex of
$\widetilde{\Gamma}$ is $\widetilde{d}(q) = 2\,d(q)$, which is
even.  By Euler's theorem, there is an Euler circuit $\Gamma'$ of
$\widetilde{\Gamma}$, which may be thought of as a closed path
which traverses each edge of $\Gamma$ exactly {\em twice}.  Now at
each of the points $\{q_1, \dots, q_d\}$ along $\Gamma'$ which are
mapped to $q \in \Gamma$, we may consider the exterior angle 
${\rm c}(q_i)$.  The sum of these exterior angles, however,
depends on the choice of the Euler circuit $\Gamma'$.  For
example, if $\Gamma$ is the union of the $x$-axis and the $y$-axis
in Euclidean space $\Re^3$, then one might choose $\Gamma'$ to
have four right angles, or to have four straight angles, or
something in between, with completely different values of total
curvature.  In order to form a version of total curvature at a
vertex $q$ which only depends on the original graph $\Gamma$ and
not on the choice of Euler circuit $\Gamma'$, it is necessary to
consider some of the exterior angles as partially balancing
others.  In the example just considered, where $\Gamma$ is the
union of two orthogonal lines, two opposing right angles will be
considered to balance each other completely, so that 
${\rm ntc}(q)=0$, regardless of the choice of Euler circuit of the
double.

It will become apparent 
that the connected character
of an Euler circuit of $\widetilde\Gamma$ is not required for what
follows.  Instead, we shall refer to a {\em parameterization}
$\Gamma'$ of the double $\widetilde\Gamma$, which is a mapping
from a $1$-dimensional manifold without boundary, not necessarily
connected;  the mapping is assumed to cover each edge of
$\widetilde\Gamma$ once.

The nature of ${\rm ntc}(q)$ is clearer when it is localized on
$S^2$, analogously to
\cite{M}.  In the case $d(q)=2$, Milnor observed that the
exterior angle at the vertex $q$ equals half
the area of
those $e \in S^2$ such that the linear function
$\langle e, \cdot \rangle$, restricted to $\Gamma$, has a local
maximum at $q$.  In our context, we may describe ${\rm ntc}(q)$ 
as one-half the integral over the sphere of the number of
{\em net local maxima}, which is half the difference of local
maxima and local minima.  Along the parameterization $\Gamma'$ of
the double of $\Gamma$, the linear function $\langle e, \cdot
\rangle$ may have a local maximum at some of the vertices $q_1,
\dots, q_d$ over $q$, and may have a local minimum at others.  In
our construction, each local minimum balances against one local
maximum.  If there are more local minima than local maxima, the
number ${\rm nlm}(e,q)$, the net number of local maxima, will be
negative;  however, our definition uses only the positive part
$[{\rm nlm}(e,q)]^+$.

       We need to show that
$$ \int_{S^2} [{\rm nlm}(e,q)]^+ \,dA_{S^2}(e)$$
is independent of the choice of parameterization, and in fact is
equal to $2 \, {\rm ntc}(q)$;  this will follow from another way
of computing ${\rm nlm}(e,q)$ (see Corollary \ref{cor2} below).

%
\begin{defi}\label{defnlm}
Let a parameterization \,$\Gamma'$ of the double of \, $\Gamma$ be
given.  Then a vertex $q$ of \,$\Gamma$ corresponds to a number of
vertices $q_1, \dots, q_d$ of \,$\Gamma'$, where $d$ is the degree 
$d(q)$ of $q$ as a vertex of \,$\Gamma$.  Choose $e \in S^2$. If 
$q \in \Gamma$ is a local extremum of $\langle e, \cdot \rangle$,
then we consider $q$
as a vertex of degree $d(q) = 2$.  Let ${\rm lmax}(e,q)$ be the
number of local maxima of $\langle e, \cdot \rangle$ along
\,$\Gamma'$ at the points $q_1, \dots, q_d$ over $q$,  and
similarly let ${\rm lmin}(e,q)$ be the
number of local minima.  We define the number of
{\em net local maxima} of $\langle e, \cdot \rangle$ at $q$ to be
$${\rm nlm}(e,q) = \frac12[{\rm lmax}(e,q) - {\rm lmin}(e,q)]$$.
\end{defi}

%
\begin{rem}
The definition of ${\rm nlm}(e,q)$ appears to depend not only on
$\Gamma$ but on a choice of the parameterization $\Gamma'$ of the
double of $\,\Gamma$:  ${\rm lmax}(e,q)$ and ${\rm lmin}(e,q)$ may
depend on the choice of $\Gamma'$.  However, we shall see in
Corollary \ref{cor1} below that the number of {\bf net} local
maxima ${\rm nlm}(e,q)$ is in fact independent of $\,\Gamma'$.
\end{rem}

%
\begin{rem}
We have included the factor $\frac12$ in the definition of
${\rm nlm}(e,q)$ in order to agree with the difference of the
numbers of local maxima and minima along a parameterization 
of $\,\Gamma$ itself,  
if $d(q)$ is even. 
\end{rem}

We shall {\bf assume} for the rest of this section that a unit vector
$e$ has been chosen, and that the linear ``height" function
$\langle e, \cdot \rangle$ has only a finite number of critical
points along $\Gamma$;  this excludes $e$ belonging to a subset of
$S^2$ of measure zero.  We shall also assume that the graph
$\Gamma$ is subdivided to include among the vertices all critical
points of the linear function $\langle e, \cdot \rangle$, with 
degree $d(q) = 2$ if $q$ is an interior point of one of
the topological edges of $\Gamma$.

%
\begin{defi}\label{updown}
Choose a unit vector $e$.  At a point $q \in \Gamma$ of degree 
$d = d(q)$, let the {\em up-degree} $d^+ = d^+(e,q)$ be the
number of edges of \,$\Gamma$ with endpoint $q$ on which
$\langle e, \cdot \rangle$ is greater (``higher") than
$\langle e, q \rangle$, the ``height" of $q$.
Similarly, let the {\em down-degree} $d^-(e,q)$ be the number of
edges along which $\langle e, \cdot \rangle$ is less than its
value at $q$.  Note that $d(q) = d^+(e,q) + d^-(e,q)$, for almost
all
$e$ in $S^2$.
\end{defi}

%
\begin{lem}\label{combin}
{\bf (Combinatorial Lemma)} For all $q \in \Gamma$ and for 
a.a.  $e\in S^2$, 	
${\rm nlm}(e,q) = \frac12[d^-(e,q) - d^+(e,q)]$.
\end{lem}
\noindent
\pf
Let a parameterization $\Gamma'$ of the double of $\Gamma$ be
chosen, with respect to which ${\rm lmax}(e,q)$ and ${\rm
lmin}(e,q)$ are defined.  Recall the assumption above, that
$\Gamma$ has been subdivided so that along each edge, the linear
function $\langle e, \cdot \rangle$ is strictly monotone.

Consider a vertex $q$ of $\Gamma$, of degree $d=d(q)$.  Then
$\Gamma'$  
has $2d$ edges with an endpoint among the points $q_1, \dots, q_d$ 
which are mapped to $q \in \Gamma$.  
On $2d^+$, resp. $2d^-$ of these edges, $\langle e,
\cdot \rangle$ is greater resp. less than $\langle e, q \rangle$.
But for each $1\leq i\leq d$, the parameterization $\Gamma'$ has
exactly two edges which meet at $q_i$.  Depending on the up/down
character of the two edges of $\Gamma'$ which meet at $q_i$,
$1\leq i\leq d$, we can count:\\ 
(+) If $\langle e, \cdot \rangle$ is greater than $\langle e, q
\rangle$ on both edges, then $q_i$ is a local minimum point;
there are ${\rm lmin}(e,q)$ of these among $q_1, \dots, q_d$.  \\ 
(-) If $\langle e, \cdot \rangle$ is less than $\langle e, q
\rangle$ on both edges, then $q_i$ is a local maximum point;
there are ${\rm lmax}(e,q)$ of these.  \\
(0) In all remaining cases, the linear function $\langle e, \cdot
\rangle$ is greater than $\langle e, q \rangle$ along one edge and
less along the other, in which case $q_i$ is not counted in
computing ${\rm lmax}(e,q)$ nor ${\rm lmax}(e,q)$; there are
$d(q)-{\rm lmax}(e,q)-{\rm lmin}(e,q)$ of these.

        Now count the individual edges of $\Gamma'$:  \\
(+) There are ${\ \rm lmin}(e,q)$ pairs of edges, each of which is
part of a local minimum, both of which are counted among the
$2 d^+(e,q)$ edges of $\Gamma'$ with
$\langle e, \cdot \rangle$ greater than $\langle e, q \rangle$.\\
(-) There are ${\ \rm lmax}(e,q)$ pairs of edges, each of which is
part of a local
maximum;  these are counted among the number $2d^-(e,q)$ of edges
of $\Gamma'$ with
$\langle e, \cdot \rangle$ less than $\langle e, q \rangle$.
Finally,\\
(0) there are $d(q)-{\rm lmax}(e,q)-{\rm lmin}(e,q)$ edges of
$\Gamma'$ which are not part of a local
maximum or minimum, with $\langle e, \cdot \rangle$
greater than $\langle e, q \rangle$;  and an equal number of edges
with $\langle e, \cdot \rangle$ less than $\langle e, q \rangle$.

Thus, the total number of these edges of $\Gamma'$ with
$\langle e, \cdot \rangle$ greater than $\langle e, q \rangle$ is
$$
2d^+=
2{\ \rm lmin}+(d-{\rm lmax}-{\rm lmin})=d+{\rm lmin}-{\rm lmax}.
$$
Similarly,
$$
2d^-=
2{\ \rm lmax}+(d-{\rm lmax}-{\rm lmin})=d+{\rm lmax}-{\rm lmin}.
$$
Subtracting gives the conclusion:
$$
{\rm nlm}(e,q):=
\frac{{\rm lmax}(e,q)-{\rm lmin}(e,q)}{2}=
\frac{d^-(e,q)-d^+(e,q)}{2}.
$$
\qed

%
\begin{cor}\label{cor1}
The number of net local maxima ${\rm nlm}(e,q)$ is independent of
the choice of parameterization $\Gamma'$ of the double of
$\Gamma$.
\end{cor}
\pf
Given a direction $e\in S^2$, the up-degree and down-degree
$d^\pm(e,q)$ at a vertex $q\in \Gamma$ are defined
independently of the choice of $\Gamma'$.
\qed

%
\begin{cor}\label{cor2}
For any $q \in \Gamma$, we have
${\rm ntc}(q) =
\frac12\int_{S^2} \Big[{\rm nlm}(e,q)\Big]^+ \,dA_{S^2}.$
\end{cor}
\pf
Consider $e \in S^2$.  In the definition \eqref{defnc} of 
${\rm ntc}(q),$ $\chi_i(e) = \pm 1$ whenever 
$\pm \langle e, T_i \rangle < 0$.  But the number of 
$1\leq i \leq d$ with $\pm \langle e, T_i \rangle < 0$ equals
$d^{\mp}(e,q)$, so that
$$\sum_{i=1}^d\chi_i(e)=d^-(e,q)-d^+(e,q)=2\,{\rm nlm}(e,q)$$
by Lemma \ref{combin}, for almost all $e \in S^2$.
\qed

%
\begin{defi}\label{defmu}
For a graph $\Gamma$ in $\Re^3$ and $e \in S^2$, define the
{\em multiplicity at $e$} as
$$ \mu(e) = \mu_\Gamma(e) = \sum\{{\rm nlm}^+(e,q): q
{\rm \ a\ vertex\ of\ } \Gamma
{\rm \ or\ a\ critical\ point\ of\ } \langle e,\cdot \rangle\}.$$
\end{defi}

Note that $\mu(e)$ is a half-integer.  Note also that in the case
when $\Gamma$ 
is a knot, or
equivalently, when $d(q) \equiv 2$, $\mu(e)$ is exactly
the integer $\mu(\Gamma, e)$, the number of local maxima of
$\langle e, \cdot \rangle$ along $\Gamma$ as defined in \cite{M},
p. 252.  

%
\begin{cor}\label{mucompare}
For almost all $e\in S^2$ and for any parameterization $\Gamma'$
of the double of $\Gamma$, 
$\mu_\Gamma(e) \leq \frac12\mu_{\Gamma'}(e).$
\end{cor}
\pf
We have $\mu_\Gamma(e) = 
\frac12\sum_q[{\rm lmax}_{\Gamma'}(e,q)-{\rm lmin}_{\Gamma'}(e,q)], 
\leq \frac12\sum_q {\rm lmax}_{\Gamma'}(e,q) =\frac12\mu_{\Gamma'}.$
\qed\\

If, in place of the positive part, we sum ${\rm nlm}(e,q)$ itself over
$q$ located above a plane orthogonal to $e$, we find a useful
quantity:

%
\begin{cor}\label{fibcard}
For almost all $s_0\in \Re$ and almost all $e \in S^2$, 
$$ 2 \sum\{{\rm  nlm}(e,q): \langle e,q \rangle > s_0\} =
\#(e,s_0),$$
the cardinality of the fiber 
$\{p\in \Gamma: \langle e,p \rangle = s_0\}$.
\end{cor}
\pf 
If $s_0 > \max_{p\in \Gamma} \langle e,p \rangle$, then
$\#(e,s_0)=0$. Now proceed downward, using Lemma \ref{combin} by
induction.
\qed\\

Note that the fiber cardinality of Corollary \ref{fibcard} is also
the value obtained for knots, where the more general ${\rm  nlm}$
may be replaced by the number of local maxima \cite{M}.

%
\begin{rem}\label{hidim}
In analogy with Corollary \ref{fibcard}, we expect that an
appropriate generalization of $\nc$ to curved polyhedral complexes
of dimension $\geq 2$ will in the future allow computation of the
homology of level sets and sub-level sets of a (generalized) Morse
function in terms of a generalization of ${\rm  nlm}(e,q)$.
\end{rem}

%
\begin{cor}\label{absnlm}
The multiplicity of a graph in direction $e\in S^2$ may also be
computed as $ \mu(e) = \frac12\sum_{q\in\Gamma}|{\rm nlm}(e,q)|$.
\end{cor}
\pf 
It follows from Corollary \ref{fibcard} with
$s_0<\min_{\Gamma}\langle e,\cdot \rangle$ that
$\sum_{q\in\Gamma} {\rm  nlm}(e,q)=0$, which is the difference of
positive and negative parts. The sum of these parts is
$\sum_{q\in\Gamma}|{\rm nlm}(e,q)|=2\mu(e).$
\qed
\vspace{1em}

It was shown in Theorem 3.1 of \cite{M}  that, in the case of
knots, 
${\mathcal C}(\Gamma)=\frac12 \int_{S^2}\mu(e)\, dA_{S^2}$, 
where Milnor refers to Crofton's formula.
We may now extend this result to {\em graphs}:

%
\begin{thm}\label{muthm}
For a (piecewise $C^2$) graph $\Gamma$ mapped into $\Re^3,$ the 
net total curvature has the following representation:
$$ \nc(\Gamma) = \frac12 \int_{S^2} \mu(e) \,dA_{S^2}(e). $$
\end{thm}
\pf
We have $ \nc(\Gamma) =
\sum_{j=1}^N {\rm ntc}(q_j) + 
\int_{\Gamma_{\rm reg}} |\vec{k}| \, ds,$
where $q_1, \dots, q_N$ are the vertices of \, $\Gamma$, including
local extrema as vertices of degree $d(q_j) = 2$, and where
$\rm{ntc}(q):=
\frac14\int_{S^2}\left[\sum_{i=1}^d\chi_i(e)\right]^+\,dA_{S^2}(e)$
by the definition (\ref{defnc}) of ${\rm ntc}(q)$.  Applying Milnor's
result to each $C^2$ edge, we have
${\mathcal C}(\Gamma_{\rm reg}) =
\frac12 \int_{S^2} \mu_{\Gamma_{\rm reg}}(e) \, dA_{S^2}$.  But
$\mu_{\Gamma}(e) = \mu_{\Gamma_{\rm reg}}(e) +
\sum_{j=1}^N \rm{nlm}^+(e,q_j)$, and the theorem follows.
\qed\\

%
\begin{cor}\label{notemb}
If $f:\Gamma \to \Re^3$ is piecewise $C^2$ but is not an
embedding, then the net total curvature $\nc(\Gamma)$ is well
defined, using the right-hand side of the conclusion of Theorem
\ref{muthm}. 
Moreover, $\nc(\Gamma)$ has the same value when points of  
self-intersection of $\Gamma$ are redefined as vertices.
\end{cor}

For $e \in S^2$, we shall use the notation $p_e:\Re^3 \to e\Re$ 
for the orthogonal projection $\langle e,\cdot \rangle$. We shall
sometimes identify $\Re$ with the one-dimensional subspace $e\Re$ 
of $\Re^3$.

\begin{cor}\label{1dsuff}
For any homeomorphism type $\{\Gamma\}$ of graphs, the infimum 
$\nc(\{\Gamma\})$ 
of net total curvature among mappings $f:\Gamma \to \Re^n$ is
assumed by a mapping $f_0:\Gamma \to \Re$.

For any isotopy class $[\Gamma]$ of embeddings 
$f:\Gamma \to \Re^3$, the infimum $\nc([\Gamma])$ of net total
curvature is assumed by a mapping $f_0:\Gamma \to \Re$ in the
closure of the given isotopy class.

Conversely, if $f_0:\Gamma \to \Re$ is in the
closure of a given isotopy class $[\Gamma]$ of embeddings into
$\Re^3$, then for all $\delta >0$ there is an embedding 
$f:\Gamma \to \Re^3$ in that isotopy class with 
$\nc(f)\leq\nc(f_0)+\delta$.
\end{cor}
\pf
Let $f:\Gamma \to \Re^3$ be any piecewise smooth mapping.
By Corollary \ref{notemb} and Corollary \ref{fibcard}, the net 
total curvature of the projection 
$p_e\circ f:\Gamma \to \Re$ of $f$ onto the line in the direction
of almost any $e\in S^2$ is given by 
$2\pi \mu(e)=\pi(\mu(e)+\mu(-e)).$ It follows from Theorem
\ref{muthm} that $\nc(\Gamma)$ is the average of $2\pi\mu(e)$ over
$e$ in $S^2$. But the half-integer-valued function $\mu(e)$ is
lower semi-continuous almost everywhere, as may be seen using
Definition \ref{defnlm}.  Let $e_0 \in S^2$ be a point where $\mu$
attains its essential infimum.  Then 
$\nc(\Gamma) \geq \pi\mu(e_0)=\nc(p_{e_0}\circ f).$
But $(p_{e_0}\circ f)e_0$ is the limit as $\varepsilon \to 0$ of
the map $f_\varepsilon$ whose projection in the direction $e_0$ is
the same as that of $f$ and is multiplied by $\varepsilon$ in all
orthogonal directions. Since $f_\varepsilon$ is isotopic to $f$,
$(p_{e_0}\circ f)e_0$ is in the closure of the isotopy class of
$f$.

Conversely, given $f_0:\Gamma \to \Re$ in the closure of a given
isotopy class, let $f$ be an embedding in that isotopy class 
uniformly close to $f_0\,  e_0$; $f_\varepsilon$ as constructed 
above converges uniformly to $f_0$
as $\varepsilon \to 0$, and 
$\nc(f_\varepsilon)\to \nc(f_0)$.
\qed

\vspace{1em}

%
\begin{defi}\label{defflat}
We call a mapping $f:\Gamma\to\Re^n$ {\em flat} (or $\nc$-flat) if 
$\nc(f)= \nc(\{\Gamma\})$, the minimum value for the topological type
of $\Gamma$, among all ambient dimensions $n$.
\end{defi}

In particular, Corollary \ref{1dsuff} above shows that for any
$\Gamma$, there is a flat mapping $f:\Gamma\to\Re$.

%
\begin{pro}\label{minmonot}
Consider a piecewise $C^2$ mapping $f_1:\Gamma \to \Re$.
There is a mapping $f_0:\Gamma \to \Re$ which 
is monotonic along the topological edges of $\Gamma$, 
has values at topological vertices of $\, \Gamma$ arbitrarily
close to those of $f_1$, and has $\nc(f_0) \leq \nc(f_1).$
\end{pro}
\pf
Any piecewise $C^2$ mapping $f_1:\Gamma \to \Re$ may be
approximated uniformly by mappings with a finite set of local
extreme
points, using the compactness of $\Gamma$.  Thus, we may assume
without loss of generality that $f_1$ has only finitely many local
extreme points.  Note that for a mapping $f:\Gamma \to \Re=\Re e$,
${\rm NTC}(f)=2\pi \mu(e)$:  hence, we only need to compare
$\mu_{f_0}(e)$ with $\mu_{f_1}(e)$.  

If $f_1$ is not monotonic on a topological edge $E$, then it has a
local extremum at a point $z$ in the interior of $E$. For
concreteness, we shall assume $z$ is a local maximum point; the
case of a local minimum is similar. Write $v, w$ for the endpoints
of $E$.  Let  $v_1$ be the closest local minimum point to $z$ on
the interval of $E$ from $z$ to $v$ (or $v_1=v$ if there is no
local minimum point between), and let $w_1$ be the closest local
minimum point to $z$ on the interval from $z$ to $w$ (or $w_1=w$).
Let $E_1\subset E$ denote the interval between $v_1$ and $w_1$.
Then $E_1$ is an interval of a topological edge of $\Gamma$,
having end points $v_1$ and $w_1$ and containing an interior point
$z$, such that $f_1$ is monotone increasing on the interval from
$v_1$ to $z$, and monotone decreasing on the interval from $z$ to
$w_1$.  By switching $v_1$ and $w_1$ if needed, we may assume that
$f_1(v_1) < f_1(w_1) < f_1(z)$. 

Let $f_0$ be equal to $f_1$ except on the interior of the interval
$E_1$, and map $E_1$ monotonically to the interval of $\Re$
between $f_1(v_1)$ and $f_1(w_1)$. Then for $f_1(w_1) < s <
f_1(z)$, the cardinality $\#(e,s)_{f_0} = \#(e,s)_{f_1} -2$. For
$s$ in all other intervals of $\Re$, this cardinality is
unchanged.  Therefore, ${\rm nlm}_{f_1}(w_1) = {\rm
nlm}_{f_0}(w_1)-1$, by Lemma \ref{combin}. This implies that ${\rm
nlm}^+_{f_1}(w_1) \geq {\rm nlm}^+_{f_0}(w_1)-1$.  Meanwhile,
${\rm nlm}_{f_1}(z)=1$, a term which does not appear in the
formula for $\mu_{f_0}$ (see Definition \ref{defmu}).
Thus $\mu_{f_0} \leq \mu_{f_1},$ and $\nc(f_0)\leq\nc(f_1)$.

Proceeding inductively, we remove each local extremum in the
interior of any edge of $\Gamma$, without increasing $\nc$.
\qed

\vspace{1em}

\section{REPRESENTATION FORMULA FOR NOWHERE-SMOOTH
GRAPHS}\label{nonsmooth}

Recall, while defining the total curvature for continuous
graphs in section \ref{deftc} above, we needed the monotonicity of
$\nc(P)$ under refinement of {\em polygonal} graphs $P$. We are
now ready to prove this.

%
\begin{pro}\label{monotmu}
Let $P$ and $\widetilde{P}$ be polygonal graphs in $\Re^3$, having
the same topological vertices, and homeomorphic to each other.
Suppose that every vertex of $P$ is also a vertex of
$\widetilde{P}$: $\widetilde{P}$ is a {\em refinement} of $P$.
Then for almost all $e \in S^2$, the multiplicity 
$\mu_{\widetilde{P}}(e) \geq \mu_P(e).$  As a consequence, 
$\nc(\widetilde{P}) \geq \nc(P)$.
\end{pro}
\pf
We may assume, as an induction step, that $\widetilde{P}$ 
is obtained from $P$ by replacing
the edge having endpoints $q_0$, $q_2$ with two edges, one having
endpoints $q_0$, $q_1$ and the other having endpoints $q_1$,
$q_2$. 
Choose $e \in S^2$.  We consider various  
cases:

If the new vertex $q_1$ satisfies 
$\langle e, q_0\rangle < \langle e, q_1\rangle < 
\langle e, q_2\rangle$, then 
${\rm nlm}_{\widetilde{P}}(e,q_i)={\rm nlm}_P(e,q_i)$ for 
$i = 0,2$ and ${\rm nlm}_{\widetilde{P}}(e,q_1)=0$, hence
$\mu_{\widetilde{P}}(e) = \mu_P(e)$. 

If $\langle e, q_0\rangle < \langle e, q_2\rangle < 
\langle e, q_1\rangle$, then 
${\rm nlm}_{\widetilde{P}}(e,q_0)={\rm nlm}_P(e,q_0)$ and
${\rm nlm}_{\widetilde{P}}(e,q_1)=1$.  The vertex $q_2$ 
requires more careful counting:  the up- and down-degree
$d_{\widetilde{P}}^\pm(e,q_2)=d_P^\pm(e,q_2) \pm 1$, so that by
Lemma \ref{combin}, 
${\rm nlm}_{\widetilde{P}}(e,q_2)={\rm nlm}_P(e,q_2)-1$.
Meanwhile, for each of the polygonal graphs, $\mu(e)$ is the sum
over $q$ of ${\rm nlm}^+(e,q)$, so the change from $\mu_P(e)$ to
$\mu_{\widetilde{P}}(e)$ depends on the value of 
${\rm nlm}_P(e,q_2)$:\\
(a) if ${\rm nlm}_P(e,q_2)\leq 0$, then 
${\rm nlm}_{\widetilde{P}}^+(e,q_2)={\rm nlm}_P^+(e,q_2)=0$;\\
(b) if ${\rm nlm}_P(e,q_2) = \frac12,$ then 
${\rm nlm}_{\widetilde{P}}^+(e,q_2)=
{\rm nlm}_P^+(e,q_2)-\frac12$;\\
(c) if ${\rm nlm}_P(e,q_2)\geq 1$, then 
${\rm nlm}_{\widetilde{P}}^+(e,q_2)=
{\rm nlm}_P^+(e,q_2)-1$.\\
Since the new vertex $q_1$ does not appear in $P$, recalling
that ${\rm nlm}_{\widetilde{P}}(e,q_1)=1$, we have 
$\mu_{\widetilde{P}}(e) - \mu_P(e) = +1, +\frac12$ or $0$ in the
respective cases (a), (b) or (c).  In any case, 
$\mu_{\widetilde{P}}(e) \geq \mu_P(e)$. 

The reverse inequality
$\langle e, q_1\rangle < \langle e, q_2\rangle < 
\langle e, q_0\rangle$ 
may be reduced to the case just above by replacing $e \in S^2$
with $-e$, since $\mu_P(-e)=-\mu_P(e)$ for any polhedral graph
$P$.
Then, depending whether ${\rm nlm}_P(e,q_2)$ is $\leq -1$,
$=-\frac12$ or $\geq 0$, we find that 
$\mu_{\widetilde{P}}(e)-\mu_P(e)= 
{\rm nlm}^+_{\widetilde{P}}(e, q_2)-{\rm nlm}^+_P(e,q_2) =
0$, $\frac12$, or $1$.  In any case,
$\mu_{\widetilde{P}}(e)\geq \mu_P(e)$.

These arguments are unchanged if $q_0$ is switched with
$q_2$.  This covers all cases except those in which equality
occurs between $\langle e,q_i \rangle$ and $\langle e,q_j \rangle$
($i\neq j$).  The set of such unit vectors $e$ form a set of
measure zero in $S^2$.  The conclusion 
$\nc(\widetilde{P}) \geq \nc(P)$ now follows from Theorem
\ref{muthm}.  \qed\\

We remark here that this step of proving the monotonicity for the
nowhere-smooth case differs from Milnor's argument for the knot
total curvature, where it was shown by two applications of 
the triangle inequality for spherical triangles.
 
Milnor extended his results for piecewise smooth knots to
continuous knots in \cite{M};  we shall carry out an analogous
extension to continuous graphs.

%
\begin{defi}\label{critpt}
We say a point $q\in \Gamma$ is {\em critical} 
relative to $e \in S^2$ when 
$q$ is a topological vertex of \,$\Gamma$ or when
$\langle e, \cdot\rangle$ is not monotone in any open 
interval of \, $\Gamma$ containing $q$.
\end{defi}

Note that at some points of a differentiable curve,
$\langle e, \cdot\rangle$ may have derivative zero but still
not be considered a critical point relative to $e$ by our 
definition.  This is appropriate to the $C^0$ category.
For a continuous graph $\Gamma$, when 
$\nc (\Gamma)$ is finite, we shall 
show that the number of critical points is finite for almost all
$e$ in $S^2$ (see Lemma \ref{fincrit} below).   

%
\begin{lem}\label{monotcvge}
Let $\Gamma$ be a continuous, finite graph in $\Re^3$, and choose
a sequence $\widehat{P_k}$ of \, $\Gamma$-approximating polygonal 
graphs with $\nc(\Gamma)= \lim_{k \rightarrow \infty}
\nc(\widehat{P_k}).$ 
Then for each
$e \in S^2$, there is a refinement $P_k$ of $\widehat{P_k}$ 
such that $\lim_{k \rightarrow \infty}\mu_{P_k}(e)$ exists in
$[0,\infty]$. 
\end{lem}
\pf 
First, for each $k$ in sequence, we refine 
$\widehat{P_k}$ to include all vertices of
$\widehat{P_{k-1}}$.  Then for all $e\in S^2$, 
$\mu_{\widehat{P_k}}(e) \geq \mu_{\widehat{P_{k-1}}}(e)$,
by Proposition \ref{monotmu}.
Second, we refine $\widehat{P_k}$ so that the arc of
$\Gamma$ corresponding to each edge of $\widehat{P_k}$ has
diameter $\leq 1/k$.  
Third, given a particular 
$e \in S^2$, for each edge $\widehat{E_k}$ of $\widehat{P_k}$, 
we add $0,1$ or $2$ points from $\Gamma$ as vertices of 
$\widehat{P_k}$ so that 
$\max_{\widehat{E_k}}\langle e,\cdot\rangle=
\max_E\langle e,\cdot\rangle$
where $E$ is the closed arc of \, $\Gamma$ corresponding to
$\widehat{E_k}$;  
and similarly so that 
$\min_{\widehat{E_k}} \langle e,\cdot\rangle=
\min_E\langle e,\cdot\rangle$.
Write $P_k$ for the result of this three-step refinement.  
Note that all vertices of $P_{k-1}$ appear among the vertices
of $P_k$.  Then by Proposition \ref{monotmu}, 
$$ \nc(\widehat{P_k}) \leq \nc(P_k) \leq \nc(\Gamma), $$
so we still have 
$\nc(\Gamma)= \lim_{k \rightarrow \infty} \nc(P_k).$

Now compare the values of
$\mu_{P_k}(e)=\sum_{q\in P_k} {\rm nlm_{P_k}}^+(e,q)$ with
the same sum for $P_{k-1}$.
Since $P_k$ is a refinement of $P_{k-1}$, we have
$\mu_{P_k}(e) \geq \mu_{P_{k-1}}(e)$ by Proposition
\ref{monotmu}. 

Therefore the values $\mu_{P_k}(e)$ are non-decreasing in $k$,
which
implies they are either convergent or properly divergent;  in the
latter case we write 
$\lim_{k \rightarrow \infty}\mu_{P_k}(e)= \infty$.
\qed

%
\begin{defi}
For a continuous graph $\Gamma$, define the {\em multiplicity} at
$e\in S^2$ as $\mu_\Gamma(e):= 
\lim_{k \rightarrow \infty}\mu_{P_k}(e) \in [0,\infty]$, 
where $P_k$ is a sequence of  $\Gamma$-approximating 
polygonal graphs, refined with respect to $e$, as given 
in Lemma \ref{monotcvge}.
\end{defi}

%
\begin{rem}
Note  that any two $\Gamma$-approximating polygonal graphs 
have a common refinement.  Hence, from the proof of Lemma 
\ref{monotcvge}, any two choices of sequences 
$\{\widehat{P_k}\}$ of $\Gamma$-approximating
polygonal graphs lead to the same value $\mu_\Gamma(e)$.
\end{rem}

%
\begin{lem}\label{a.a.e}
Let $\Gamma$ be a continuous, finite graph in $\Re^3$. 
Then $\mu_\Gamma: S^2 \to [0,\infty]$ takes its values in the
half-integers, or $+ \infty$.
Now assume $\nc(\Gamma) < \infty$.  Then $\mu_\Gamma$ is
integrable, hence finite almost everywhere on $S^2$, and 
\begin{equation}\label{nc=int}
 \nc(\Gamma) = \frac12 \int_{S^2} \mu_\Gamma(e) \, dA_{S^2}(e).
\end{equation}
For almost all $e \in S^2$, a sequence $P_k$ of
\, $\Gamma$-approximating polygonal graphs, converging uniformly
to $\Gamma$, may be chosen (depending on $e$) so that each
local extreme point $q$ of $\langle e,\cdot\rangle$ along $\Gamma$ 
occurs as a vertex of $P_k$ for sufficiently large $k$.
\end{lem}
\pf
Given $e\in S^2$, let $\{P_k\}$ be the sequence of	
$\Gamma$-approximating polygonal graphs from Lemma	
\ref{monotcvge}.  If $\mu_\Gamma(e)$ is finite, then	
$\mu_{P_k}(e)=\mu_\Gamma(e)$ for $k$ sufficiently 	
large, a half-integer.					

Suppose $\nc(\Gamma) < \infty$. Then			
the half-integer-valued functions $\mu_{P_k}$ are non-negative,
integrable on $S^2$ with bounded integrals since 
$\nc(P_k) < \nc(\Gamma) < \infty$, and monotone increasing in $k$.  
Thus for almost all $e \in S^2$,
$\mu_{P_k}(e) = \mu_\Gamma(e)$ for $k$ sufficiently large.

Since the functions $\mu_{P_k}$ are non-negative and pointwise
non-decreasing almost everywhere on $S^2$, it now follows 
from the Monotone Convergence Theorem that 
$$\int_{S^2} \mu_\Gamma(e)\, dA_{S^2}(e) = 
\lim_{k\to \infty}\int_{S^2} \mu_{P_k}(e)\, dA_{S^2}(e)=
2 \nc(\Gamma).$$

Finally, 
the polygonal graphs $P_k$ 
have maximum edge length $\to 0$.  
For almost all $e \in S^2$, $\langle e,\cdot\rangle$ is
not constant along any open arc of  $\Gamma$, and 
$\mu_\Gamma(e)$ is finite.  Given such an $e$,
choose $\ell = \ell(e)$ sufficiently 
large that $\mu_{P_k}(e) = \mu_\Gamma(e)$ and 
$\mu_{P_k}(-e) = \mu_\Gamma(-e)$ for all $k \geq \ell$.  Then for
$k \geq \ell$, along any edge $E_k$ of $P_k$ with corresponding
arc $E$ of $\Gamma$, the maximum and
minimum values of $\langle e,\cdot\rangle$ along $E$ occur at the
endpoints, which are also the endpoints of $E_k$.  Otherwise, as
$P_k$ is further refined, new interior local maximum resp. local
minimium points of $E$ would
contribute a new, positive value to $\mu_{P_k}(e)$ resp. to
$\mu_{P_k}(-e)$ as $k$
increases.  Since the diameter of the corresponding arc $E$ of
$\Gamma$ tends to zero as $k \to \infty$, any local
maximum or local minimum of $\langle e,\cdot\rangle$ must become
an endpoint of some edge of $P_k$ for $k$ sufficiently large, and
for $k\geq \ell$ in particular.
\qed\\

Our next lemma focuses on the regularity of a graph $\Gamma$,
originally only assumed continuous, provided it has finite net
total curvature, or another notion of total curvature of a graph 
which includes the total curvature of the edges.

%
\begin{lem}\label{fincrit}
Let $\Gamma$ be a continuous, finite graph in $\Re^3$, with
$\nc(\Gamma)<\infty$. Then $\Gamma$ has continuous one-sided unit
tangent vectors $T_1(p)$ and $T_2(p)$ at each point $p$, not a 
topological vertex.  If $p$ is a
vertex of degree $d$, then each of the $d$ edges which meet at
$p$ have well-defined unit tangent vectors at $p$:
$T_1(p),\dots,T_d(p)$.  For almost all $e \in S^2$, 
\begin{equation}\label{mu=sum}
\mu_\Gamma(e) = \sum_q\{{\rm nlm}(e,q)\}^+,
\end{equation}
where the sum is over the {\em finite} number of topological
vertices of \, $\Gamma$ and critical points $q$ of 
$\langle e, \cdot \rangle$ along $\Gamma$. 
Further, for each $q$,
${\rm nlm}(e,q)= \frac12[d^-(e,q) - d^+(e,q)]$.
All of these critical points which are not topological vertices 
are local extrema of $\langle e,\cdot\rangle$ along $\Gamma$.  
\end{lem}
\pf
We have seen in the proof of Lemma \ref{a.a.e} that for almost all 
$e \in S^2$, the linear function $\langle e,\cdot\rangle$ is 
not constant along any open arc of  $\Gamma$, and by Lemma
\ref{monotcvge}
there is a sequence $\{P_k\}$ of $\Gamma$-approximating 
polygonal graphs with
$\mu_\Gamma(e) = \mu_{P_k}(e)$ for $k$ sufficiently large.
We have further shown that each local maximum point of 
$\langle e,\cdot\rangle$ is a vertex of $P_k$, possibly of degree 
two,  for $k$ large enough.  Recall that
$\mu_{P_k}(e) = \sum_q{\rm nlm}_{P_k}^+(e,q)$.
Thus, each local maximum point $q$ for 
$\langle e, \cdot \rangle$ along $\Gamma$ provides a non-negative
term ${\rm nlm}_{P_k}^+(e,q)$ in the sum for $\mu_{P_k}(e)$.  
Fix such an integer $k$.

Consider a point $q\in \Gamma$ which is not a topological 
vertex of $\Gamma$ but is a critical point of 
$\langle e,\cdot\rangle$.  We shall show, by an argument similar 
to one used by van Rooij in \cite{vR}, that $q$ 
must be a local extreme point.
As a first step, we show that $\langle e,\cdot\rangle$ is monotone
on a sufficiently small interval on either side of $q$.
Choose an ordering of the closed edge $E$ of 
$\Gamma$ containing $q$, and consider the interval $E_+$ of 
points $\geq q$ with respect to this ordering.  Suppose that 
$\langle e,\cdot\rangle$ is not monotone on any subinterval of
$E_+$ with $q$ as endpoint.
Then in any interval $(q,r_1)$ there are points $p_2 > q_2 > r_2$
so that the numbers 
$\langle e,p_2 \rangle,\langle e,q_2\rangle,\langle e, r_2\rangle$
are not monotone.
It follows by an induction argument that there exist decreasing 
sequences $p_n \to q$, $q_n \to q$, and 
$r_n \to q$ of points of $E_+$ such that for each $n$, 
$r_{n-1} > p_n > q_n > r_n > q$, but the value 
$\langle e,q_n\rangle$ lies outside of the closed interval between 
$\langle e,p_n\rangle$ and $\langle e,r_n\rangle$.  As a
consequence, 
there is a local extremum $s_n \in (r_n, p_n)$.  Since 
$r_{n-1} > p_n$, the $s_n$ are all distinct, $1\leq n < \infty$.  
But by Lemma \ref{a.a.e},
all local extreme points, specifically $s_n$, of 
$\langle e, \cdot \rangle$ along $\Gamma$ 
occur among the {\em finite}
number of vertices of $P_k$, a contradiction.  This shows that 
$\langle e, \cdot \rangle$ is monotone on an interval to the right
of $q$.  An 
analogous argument shows that $\langle e, \cdot \rangle$
is monotone on an interval to the left of $q$.

Recall that  for a {\em critical point} $q$ relative to $e$, 
$\langle e,\cdot\rangle$ is not monotone on any neighborhood 
of $q$.  
Since $\langle e,\cdot\rangle$ is monotone on an interval on
either side, the sense of monotonicity 
must be opposite on the two sides of $q$.  Therefore every
critical point $q$ along $\Gamma$ for $\langle e, \cdot \rangle$, 
which is not a topological vertex, is a local extremum.

We have chosen $k$ large enough that $\mu_\Gamma(e) = \mu_{P_k}(e)$.
Then for
any edge $E_k$ of $P_k$, the function $\langle e, \cdot \rangle$
is monotone along the corresponding arc $E$ of $\Gamma$, as
well
as along $E_k$.  Also, $E$ and $E_k$ have common end points.
It follows that for each $t \in \Re$, the cardinality  
$\#(e,t)$ of the fiber $\{q\in \Gamma: \langle e,q \rangle =t \}$ 
is the same for $P_k$ as for $\Gamma$.  We may see from 
Lemma \ref{combin} applied to $P_k$ that for each vertex 
or critical point $q$,
${\rm nlm}_{P_k}(e,q) = \frac12[d_{P_k}^-(e,q) - d_{P_k}^+(e,q)]$;  
but ${\rm nlm}(e,q)$ and $d^\pm(e,q)$ have the {\it same} values 
for $\Gamma$ as for $P_k$.  The formula 
$\mu_\Gamma(e) = \sum_q\{{\rm nlm}_\Gamma(e,q)\}^+$ now follows 
from the  
corresponding formula for $P_k$, for almost all $e \in S^2$.

Consider an open interval $E$ of $\Gamma$ with endpoint $q$. We
have just shown that for a.a. $e\in S^2$,
$\langle e, \cdot \rangle$ is monotone on a subinterval with
endpoint $q$. Choose a sequence $p_\ell$ from $E$, $p_\ell\to q$,
and write $T_\ell := \frac{p_\ell-q}{|p_\ell-q|} \in S^2$. Then 
$\lim_{\ell\to\infty}T_\ell$ exists. Otherwise, since $S^2$ is
compact, there are subsequences $\{T_{m_n}\}$ and $\{T_{k_n}\}$
with $T_{m_n} \to T'$ and $T_{k_n} \to T'' \neq T'$. But for an
open set of $e\in S^2$, 
$\langle e, T' \rangle < 0 < \langle e, T'' \rangle$. For such $e$,
$\langle e, q_{m_n}\rangle<\langle e, q\rangle<\langle e, q_{k_n} \rangle$
for $n >> 1$. That is, as $p\to q$, $p\in E$, 
$\langle e, p \rangle$ assumes values above and below 
$\langle e, q\rangle$ infinitely often, contradicting monotonicity
on an interval starting at $q$ for a.a. $e\in S^2$.

This shows that $\Gamma$ has one-sided tangent vectors 
$T_1(q), \dots, T_d(q)$ at each point $q\in \Gamma$ of degree 
$d=d(q)$ ($d=2$ if q is not a topological vertex). Further, as
$k\to \infty,$ $T_i^{P_k}(q) \to T_i^\Gamma(q)$, 
$1\leq i\leq d(q)$, since edges of $P_k$ have diameter 
$\leq \frac{1}{k}$.

The remaining conclusions follow readily.

\qed

%
\begin{cor}\label{ctsnc}
Let $\Gamma$ be a continuous, finite graph in $\Re^3$, with
$\nc(\Gamma)<\infty$. Then for each  point $q$ of
$\Gamma$, the contribution at $q$ to net total curvature is
given by equation \eqref{defnc}, where for $e \in S^2$, 
$\chi_i(e)=$ the sign of $\langle -T_i(q), e \rangle$, 
$1\leq i \leq d(q)$.  
(Here, if $q$ is not a topological vertex, 
we understand $d=2$.)
\end{cor}
\pf
According to Lemma \ref{fincrit}, for $1\leq i \leq d(q)$,
$T_i(q)$ is defined and tangent to an edge $E_i$ of $\Gamma$,
which is continuously differentiable at its end point $q$.  
If $P_n$ is a sequence of
$\Gamma$-approximating polygonal graphs with maximum edge length
tending to $0$, then the corresponding unit tangent vectors
$T^{P_n}_i(q) \to T^{\Gamma}_i(q)$ as $n \to \infty$.  For each
$P_n$, we have 
$$
{\rm ntc}^{P_n}(q) =
\frac{1}{4}\int_{S^2}
\left[\sum_{i=1}^d{\chi_i}^{P_n}(e)\right]^+\,dA_{S^2}(e),
$$
and ${\chi_i}^{P_n} \to {\chi_i}^\Gamma$ in measure on $S^2$.
Hence, the integrals for $P_n$ 
converge to those for $\Gamma$, which is equation \eqref{defnc}.
\qed\\

We are ready to state the formula for net total curvature,
by localization on $S^2$, a generalization of Theorem \ref{muthm}:

%
%
\begin{thm}\label{muthm2}
For a continuous graph $\Gamma,$ the net total curvature 
$\nc(\Gamma) \in (0,\infty]$ has 
the following representation:
$$ \nc(\Gamma) = \frac14 \int_{S^2} \mu(e) \,dA_{S^2}(e), $$
where, for almost all $e \in S^2$, the multiplicity 
$\mu(e)$ is a positive half-integer or $+\infty$, given as the
finite sum \eqref{mu=sum}. 
\end{thm}

\pf
If $\nc(\Gamma)$ is finite, then the theorem follows from Lemma 
\ref{a.a.e} and Lemma \ref{fincrit}.

Suppose $\nc(\Gamma) = \sup \nc(P_k)$ is infinite, where $P_k$ 
is a refined sequence of polygonal graphs as in Lemma
\ref{monotcvge}.  Then $\mu_\Gamma(e)$ is the non-decreasing 
limit of $\mu_{P_k}(e)$ for all $e \in S^2$.  Thus    
$\mu_\Gamma(e) \geq \mu_{P_k}(e)$ for all $e$ and $k$, and
$\mu_\Gamma(e) = \mu_{P_k}(e)$ for $k\geq\ell(e)$.  This implies
that $\mu_\Gamma(e)$ is a positive half-integer or $\infty$.
Since $\nc(\Gamma) = \infty$, the integral 
$$ \nc(P_k) = \frac12\int_{S^2}\mu_{P_k}(e) \,dA_{S^2}(e)$$
is arbitrarily large as $k \to \infty$, but for each $k$ is 
less than or equal to 
$$ \frac12 \int_{S^2} \mu_\Gamma(e) \,dA_{S^2}(e).$$
Therefore this latter integral equals $\infty$, and thus equals
$\nc(\Gamma).$
\qed\\

We turn our attention next to the tameness of graphs of finite
total curvature.

%
\begin{pro}\label{untangle}
Let $n$ be a positive integer, and write $Z$ for the set of $n$-th
roots of unity in $\C = \Re^2$.  Given a continuous one-parameter 
family $S_t$, $ 0 \leq t < 1$, of sets of $n$ points in 
$\Re^2$, there exists a continuous
one-parameter family $\Phi_t:\Re^2 \to \Re^2$ of homeomorphisms
with compact support such that $\Phi_t(S_t) = Z$, $0 \leq t < 1$. 
\end{pro}

\pf
It is well known that there is an isotopy $\Phi_0:\Re^2 \to \Re^2$
such that $\Phi_0(S_0) = Z$ and $\Phi_0 =$ id outside of a
compact set.  
This completes the case 
$ t_0 = 0$ of the following continuous induction argument.

Suppose that $[0,t_0] \subset [0,1)$ is a subinterval such that
there exists a continuous one-parameter family
$\Phi_t:\Re^2 \to \Re^2$ of homeomorphisms with compact support,
with $\Phi_t(S_t) = Z$ for all $0 \leq t \leq t_0$.  
We shall extend this property to an interval $[0,t_0+\delta]$.
Write $B_\varepsilon (Z)$ for the union of balls		
$B_\varepsilon (\zeta_i)$ centered at the $n$ roots of unity	
$\zeta_1, \dots \zeta_n$.  For $\varepsilon <			
\sin{\frac{\pi}{n}},$ these balls are disjoint.  		
We may choose $0 < \delta < 1-t_0$  such that
$\Phi_{t_0}(S_t) \subset B_\varepsilon(Z)$ for all
$t_0 \leq t \leq t_0 + \delta.$  Write the points of $S_t$ as
$x_i(t), \ 1 \leq i \leq n,$ where
$\Phi_{t_0}(x_i(t)) \in B_\varepsilon(\zeta_i)$.  For each
$t \in [t_0, t_0 + \delta],$ each of the balls
$B_\varepsilon(\zeta_i)$ may be mapped onto itself by a 
homeomorphism $\psi_t$, varying continuously with $t$, such that
$\psi_{t_0}$ is the identity, $\psi_t$ is the identity near the
boundary of $B_\varepsilon(\zeta_i)$ for all 
$t \in [t_0, t_0 + \delta]$, and
$\psi_t(\Phi_{t_0}(x_i(t))) = \zeta_i$ for all such $t$. 
For example, we may construct $\psi_t$ so that
for each $y \in B_\varepsilon(\zeta_i)$, $y-\psi_t(y)$ is
parallel to $\Phi_{t_0}(x_i(t)) - \zeta_i$.  We now define
$\Phi_t = \psi_t \circ \Phi_{t_0}$ for each
$t \in [t_0, t_0 + \delta].$

As a consequence, we see that there is no maximal interval
$[0,t_0] \subset [0, 1)$ such that
there is a continuous one-parameter family
$\Phi_t:\Re^2 \to \Re^2$ of homeomorphisms with compact support
with $\Phi_t(S_t) = Z$, for all $0 \leq t \leq t_0$.
Thus, this property holds for the entire interval $0 \leq t<1$.
\qed\\

In the following theorem, the total curvature of a graph may be
understood in terms of any definition which includes the total
curvature of edges and which is continuous as a function of the
unit tangent vectors at each vertex.  This includes net total
curvature, TC of \cite{T} and CTC of \cite{GY1}.

%
%
\begin{thm}\label{tame}
Suppose $\Gamma \subset \Re^3$ is a continuous graph with finite
total curvature.  Then for any $\varepsilon > 0$, $\Gamma$ is 
isotopic to a $\Gamma$-approximating polygonal graph $P$ with
edges of length at most $\varepsilon$, whose total curvature is 
less than or equal to that of $\Gamma$.
\end{thm}
\pf
Since $\Gamma$ has finite total curvature, by Lemma \ref{fincrit},
at each topological vertex of degree $d$ the edges have
well-defined unit tangent vectors $T_1, \dots, T_d$, which are
each the limit of the unit tangent vectors to the corresponding
edges.
If at each vertex the unit tangent vectors $T_1, \dots, T_d$ 
are distinct, then any sufficiently fine
$\Gamma$-approximating polygonal graph will be isotopic to
$\Gamma$;  this easier case is proven.

We consider therefore $n$ edges $E_1, \dots, E_n$ which end at a
vertex $q$ with common unit tangent vectors $T_1 = \dots = T_n$.
Choose orthogonal coordinates $(x,y,z)$ for $\Re^3$ so that this 
common tangent vector $T_1=\dots=T_n=(0,0,-1)$ and $q = (0,0,1)$.
For some $\varepsilon>0,$ in the slab 			
$1-\varepsilon\leq z\leq 1,$ the edges $E_1, \dots, E_n$  
project one-to-one onto the $z$-axis. 			
After rescaling about $q$ by a 
factor $\geq \frac{1}{\varepsilon}$, $E_1, \dots, E_n$ form a 
braid $B$ of $n$ strands in the 
slab $0 \leq z < 1$ of $\Re^3$, plus the point $q=(0,0,1)$.  Each
strand $E_i$  has $q$ as an endpoint, and the coordinate $z$ is 
strictly monotone along $E_i$, $1\leq i \leq n$.  Write 
$S_t = B \cap \{ z = t\}$.
Then $S_t$ is a set of $n$ distinct points in the plane 
$\{ z = t\}$ for each $0\leq t<1$.  According to Proposition
\ref{untangle}, there are homeomorphisms $\Phi_t$ of the plane
$\{ z = t\}$ for each $0\leq t<1$, isotopic to the identity in
that plane,  continuous as a function of $t$,
such that $\Phi_t(S_t) = Z \times \{t\},$ where $Z$ is the set of
$n$th roots of unity in the $(x,y)$-plane, and $\Phi_t$ is the
identity outside of a compact set of the plane $\{ z = t\}$.  

We may suppose that $S_t$ lies in the open disk of radius $a(1-t)$
of the plane $\{ z = t\}$, for some (arbitrarily small) 
constant $a>0$. We modify
$\Phi_t$, first replacing its values with $(1-t) \Phi_t$ inside
the disk of radius $a(1-t)$.  We then modify $\Phi_t$ outside the
disk of radius $a(1-t)$, such that $\Phi_t$ is the identity
outside the disk of radius $2a(1-t)$.

Having thus modified the homeomorphisms $\Phi_t$ of the planes
$\{ z = t\}$, we may now define an isotopy $\Phi$ of $\Re^3$ by
mapping each plane $\{ z = t\}$ to itself by the homeomorphism
$\Phi_0^{-1} \circ \Phi_t$, $0\leq t<1$; and extend to the
remaining planes $\{ z = t\}$, $t\geq 1$ and $t<0$, by the 
identity.  Then the closure
of the image of the braid $B$ is the union of line segments from
$q =(0,0,1)$ to the $n$ points of $S_0$ in the plane $\{ z = 0\}$.
Since each $\Phi_t$ is isotopic to the identity
in the plane $\{ z = t\}$, $\Phi$ is isotopic to
the identity of $\Re^3$.

This procedure may be carried out in disjoint sets of $\Re^3$
surrounding each unit vector which occurs as tangent 
vector to more than one edge at a vertex of $\Gamma$.  Outside
these sets, we inscribe a polygonal arc in each edge of $\Gamma$
to obtain a $\Gamma$-approximating polygonal graph $P$.  By
Definition \ref{gendefnet}, $P$ has total curvature less than or
equal to the total curvature of $\Gamma$.
\qed

\vspace{1em}

Artin and Fox \cite{AF} introduced the notion of {\em tame} and
{\em wild} knots in $\Re^3$;  the extension to graphs is the
following

%
\begin{defi}
We say that a graph in $\Re^3$ is {\em tame} if it is isotopic 
to a polyhedral graph;  otherwise, it is {\em wild}.
\end{defi}

Milnor proved in \cite {M} that knots of finite total curvature
are tame. More generally, we have

%
\begin{cor}
A continuous graph $\Gamma \subset \Re^3$ of finite total
curvature is tame.
\end{cor}
\pf
This is an immediate consequence of Theorem \ref{tame}, since the
$\Gamma$-approximating polygonal graph $P$ is isotopic to
$\Gamma$.
\qed

%
\begin{obs}
Tameness does not imply finite total curvature.
\end{obs}
For a well-known example, consider $\Gamma\subset\Re^2$ to be
the continuous curve
$\{(x,h(x)):x\in [-1,1]\}$ where the function
$$ h(x)=-\frac{x}{\pi}\sin\frac{\pi}{x},$$
$h(0)=0$, has a sequence of zeroes $\pm\frac{1}{n}\to 0$ as
$n\to\infty$. Then
the total curvature of $\Gamma$ between $(0,\frac{1}{n})$ and
$(0,\frac{1}{n+1})$ converges to $\pi$ as $n \to\infty$. Thus 
${\mathcal C}(\Gamma) = \infty$.

On the other hand, $h(x)$ is continuous on $[-1,1]$, from which it
readily follows that $\Gamma$ is tame.

%
\section{ON VERTICES OF SMALL DEGREE}\label{three/four}

We are now in a position to illustrate some properties of net
total
curvature $\nc(\Gamma)$ in a few relatively simple cases, and to
make some observations regarding $\nc(\{\Gamma\})$, the minimum
net total curvature for the homeomorphism type of a graph 
$\Gamma \subset \Re^n$ (see Definition \ref{defflat} above).

%
\subsection{Minimum curvature for given degree}

%
\begin{pro}\label{val3}
If a vertex $q$ has {\bf odd} degree, then 
$\rm{ntc}(q) \geq \pi/2$.  If $d(q)=3$, then equality holds if and 
only if the three tangent vectors $T_1, T_2, T_3$ at $q$ are
coplanar but do not lie in any open half-plane. If $q$ has 
{\bf even} degree $2m$, then the minimum value of $\rm{ntc}(q)$
is $0$. Moreover, the equality $\rm{ntc}(q)=0$ only occurs when 
$T_1(q), \dots, T_{2m}(q)$ form $m$ opposite pairs. 
\end{pro}

\pf 
Let $q$ have odd degree $d(q)=2m+1$. Then from Lemma
\ref{combin}, 
for any $e\in S^2$, we see that ${\rm nlm}(e,q)$ is a
half-integer $\pm \frac12, \dots, \pm \frac{2m+1}{2}$. 
In particular, $|{\rm nlm}(e,q)|\geq \frac12$. 
Corollary \ref{cor2} and the proof of 
Corollary \ref{absnlm} show that
$${\rm ntc}(q) =
\frac{1}{4}\int_{S^2} \Big|{\rm nlm}(e,q)\Big|\,dA_{S^2}. $$
Therefore ${\rm ntc}(q) \geq \pi/2$. 

If the degree $d(q)=3$, then
$|{\rm nlm}(e,q)|= \frac12$ if and only if both $d^+(q)$ and
$d^-(q)$ are nonzero, that is, $q$ is not a local
extremum for $\langle e,\cdot \rangle$. 
If $\rm{ntc}(q) = \pi/2$, then this must be true for
almost every direction $e\in S^2$.
Thus, the three tangent vectors
must be coplanar, and may not lie in an open half-plane.

If $d(q)=2m$ is even and equality $\rm{ntc}(q)=0$ holds, then the
formula above for ${\rm ntc}(q)$ in terms of $|{\rm nlm}(e,q)|$ 
would require ${\rm nlm}(e,q)\equiv 0$, 
and hence $d^+(e,q)=d^-(e,q)=m$ for almost all $e\in S^2$:
whenever $e$
rotates so that the plane orthogonal to $e$ passes $T_i$,
another tangent vector $T_j$ must cross the plane in the
opposite direction, for a.a. $e$, which implies $T_j=-T_i$.
\qed

%
\begin{obs}\label{odd>3}
If a vertex $q$ of odd degree $d(q)=2p+1$, 
has the minimum value $\rm{ntc}(q)=\pi/2$,
and a hyperplane $P\subset \Re^n$ contains an even number
of the tangent vectors at $q$, and no others, then these tangent 
vectors form opposite pairs. 
\end{obs}

The proof is seen by fixing any
$(n-2)$-dimensional subspace $L$ of $P$ and rotating $P$ by a
small positive or negative angle $\delta$ to a hyperplane
$P_\delta$ containing $L$. 
Since $P_\delta$ must have $k$ of the vectors 
$T_1, \dots, T_{2p+1}$ on one side and $k+1$ on the other side,
for some $0\leq k \leq p$, by comparing $\delta>0$ with $\delta<0$
it follows that exactly half of the tangent vectors in $P$ lie
nonstrictly on each side of $L$. The proof may be continued as in
the last paragraph of the proof of Proposition \ref{val3}. In
particular, any 
two independent tangent vectors $T_i$ and $T_j$ share
the $2$-plane they span with a third, the three vectors not lying
in any open half-plane: in fact, the third vector needs to lie in
any hyperplane containing $T_i$ and $T_j$.

For example, a flat $K_{5,1}$ in $\Re^3$ must have five straight
segments, two being opposite; and the remaining three being
coplanar but not in any open half-plane. This includes the case of
four coplanar line segments, since the four must be in opposite
pairs, and either opposing 
pair may be considered as coplanar with the fifth segment.

%
\subsection{Non-monotonicity of $\nc$ for subgraphs}

%
\begin{obs}\label{notmonotone}
If $\Gamma_0$ is a subgraph of a graph $\Gamma$, then
$\nc(\Gamma_0)$ might {\bf not} be $\leq \nc(\Gamma).$ 
\end{obs}

For a simple polyhedral example, we may
consider the ``butterfly" graph $\Gamma$ in the plane with six
vertices:  $q_0^\pm = (0,\pm 1), q_1^\pm = (1,\pm 3),$ and
$q_2^\pm = (-1,\pm 3)$. $\Gamma$ has seven edges: three vertical
edges $L_0, L_1$ and $L_2$ are the line segments $L_i$ joining
$q_i^-$ to $q_i^+$. Four additional edges are the line segments
from $q_0^\pm$ to $q_1^\pm$ and from $q_0^\pm$ to $q_2^\pm$, which
form the smaller angle $2 \alpha$ at $q_0^\pm$, where $\tan \alpha
= 1/2$, so that $\alpha < \pi/4.$

The subgraph $\Gamma_0$ will be $\Gamma$ minus the interior of
$L_0$.  Then $\nc(\Gamma_0) = {\mathcal C}(\Gamma_0)= 6 \pi - 8
\alpha.$ However, 
$\nc(\Gamma) = 4(\pi - \alpha) + 2(\pi/2) = 5 \pi - 4 \alpha,$
which is $<\nc(\Gamma_0).$
\qed\\

The monotonicity property, which is shown in Observation 
\ref{notmonotone} to fail for $\nc(\Gamma)$, is a
virtue of Taniyama's total curvature $\tc(\Gamma)$.

%
\subsection{Net total curvature $\neq$ cone total curvature
$\neq$ Taniyama's total curvature}

\vspace{1em}

It is not difficult to construct three unit vectors 
$T_1, T_2, T_3$ in $\Re^3$ 
such that the values of ${\rm ntc}(q)$,
${\rm ctc}(q)$ and ${\rm tc}(q)$, with these vectors as the 
$d(q) = 3$ tangent vectors to a graph at a vertex $q$, have 
different values.  For example, we
may take $T_1, T_2$ and $T_3$ to be
three unit vectors in a plane, making equal angles $2\pi/3$.
According to Proposition \ref{val3}, we have the contribution to
net total curvature ${\rm ntc}(q) = \pi/2$.  But the contribution
to cone total curvature is ${\rm ctc}(q) = 0$. Namely,
${\rm ctc}(q) := \sup_{e \in S^2}
\sum_{i=1}^3\left(\frac{\pi}{2}-\arccos\langle T_i, e\rangle
\right).$  In this supremum, we may choose $e$ to be normal 
to the plane of $T_1, T_2$ and $T_3$, and ${\rm ctc}(q) = 0$
follows.
Meanwhile, ${\rm tc}(q)$ is the sum of the exterior angles formed
by the three pairs of vectors, each equal to $\pi/3$, so that
${\rm tc}(q) = \pi$.

A similar computation for degree $d$ and coplanar vectors making
equal angles gives ${\rm ctc}(q) = 0$, and 
${\rm tc}(q) = \frac{\pi}{2}\Big[\frac{(d-1)^2}{2}\Big]$ (brackets
denoting integer part), while ${\rm ntc}(q) = \pi/2$ for $d$ odd,
${\rm ntc}(q) = 0$ for $d$ even.  This example indicates that 
${\rm tc}(q)$ may be significantly larger than ${\rm ntc}(q)$.
In fact, we have 

%
\begin{obs}\label{tc>>ntc}
If a vertex $q$ of a graph $\Gamma$ has degree $d=d(q)\geq 2$, then
${\rm tc}(q) \geq (d-1) {\rm ntc}(q)$.
\end{obs}

This 
follows from the definition \eqref{defnc} of
${\rm ntc}(q)$.  Let $T_1, \dots, T_d$ be the unit tangent vectors
at $q$.  The exterior angle between $T_i$ and $T_j$ is 
$$\arccos\langle -T_i,T_j \rangle =
\frac{1}{4} \int_{S^2} (\chi_i + \chi_j)^+ \, dA_{S^2}.$$
The contribution ${\rm tc}(q)$ at $q$ to total curvature
$\tc(\Gamma)$ equals the sum of these integrals over all 
$1 \leq i < j \leq d$.  The sum of the integrands is 
$$ \sum_{1\leq i<j\leq d}(\chi_i + \chi_j)^+\geq
\Bigg[\sum_{1\leq i<j\leq d}(\chi_i + \chi_j)\Bigg]^+ = 
(d-1)\Big[\sum_{i=1}^d \chi_i\Big]^+.$$
Integrating over $S^2$ and dividing by $4$, we have 
${\rm tc}(q) \geq (d-1) {\rm ntc}(q)$.
\qed

\vspace{1em}

%
\subsection{Conditional additivity of net total curvature under 
taking union}

\vspace{1em}

Observation 3 shows the failure of monotonicity of
$\nc$ for
subgraphs due to the cancellation phenomena at each vertex.  The
following subadditivity statement specifies the necessary and
sufficient condition for the additivity of net total curvature
under taking union of graphs.

%
\begin{pro}\label{subadditivity}
Given two graphs $\Gamma_1$ and $\Gamma_2\subset \Re^n$ with 
$\Gamma_1 \cap \Gamma_2 = \{p_1, \dots, p_N\}$, the net total 
curvature of $\,\Gamma = \Gamma_1 \cup \Gamma_2$ obeys the 
sub-additivity law 
\begin{eqnarray}\label{subadd}
\nc(\Gamma) & = & \nc(\Gamma_1) + \nc(\Gamma_2)+\nonumber\\
&+&\frac12\sum_{j=1}^N\int_{S^2}[{\rm nlm}_{\Gamma}^+(e, p_j)-
{\rm nlm}_{\Gamma_1}^+(e, p_j)-
{\rm nlm}_{\Gamma_2}^+(e, p_j)]\, dA_{S^2} \\
& \leq & \nc(\Gamma_1) + \nc(\Gamma_2).\nonumber
\end{eqnarray}

\noindent
In particular, additivity holds if and only if
\[
{\rm nlm}_{\Gamma_1}(e, p_j)\, {\rm nlm}_{\Gamma_2}(e, p_j) \geq 0
\]
for all points $p_j$ of $\,\Gamma_1 \cap \Gamma_2$  and 
almost all $e \in S^2$.
\end{pro}
\pf
The edges of $\Gamma$ and vertices other than $p_1, \dots, p_N$
are edges and vertices of $\Gamma_1$ or of $\Gamma_2$, so we only
need to consider the contribution at the vertices 
$p_1, \dots, p_N$ to $\mu(e)$ for $e\in S^2$ 
(see Definition \ref{defmu}).  			
The sub-additivity
follows from the general inequality $(a+b)^+ \leq a^+ + b^+$ for
any real numbers $a$ and $b$.  Namely, let 
$a:= {\rm nlm}_{\Gamma_1}(e, p_j)$ and 
$b:= {\rm nlm}_{\Gamma_2}(e, p_j)$, so that 
${\rm nlm}_{\Gamma}(e, p_j) = a+b$, as follows from Lemma 
\ref{combin}.  Now integrate both sides of the inequality over
$S^2$, sum over $j=1, \dots, N$ and apply Theorem \ref{muthm}.

As for the equality case, suppose that 
$ab\geq 0$. 
We then note that either 
$a > 0 $ \& $b > 0$, or 
$a < 0 $ \& $b < 0$, or 
$a = 0 $, or $b = 0$.  In all four cases, we have 
$a^+ + b^+ = (a+b)^+$. Applied with 
$a= {\rm nlm}_{\Gamma_1}(e, p_j)$
and $b= {\rm nlm}_{\Gamma_2}(e, p_j)$, assuming that 
${\rm nlm}_{\Gamma_1}(e, p_j){\rm nlm}_{\Gamma_2}(e,p_j)\geq 0$
holds for all $j = 1, \dots, N$ and almost
all $e \in S^2$, this implies that
$\nc(\Gamma_1\cup\Gamma_2)=\nc(\Gamma_1)+\nc(\Gamma_2).$

To show that the equality 
$\nc(\Gamma_1 \cup\Gamma_2)=\nc(\Gamma_1)+\nc(\Gamma_2)$ 
implies the inequality   
${\rm nlm}_{\Gamma_1}(e, p_j){\rm nlm}_{\Gamma_2}(e, p_j)\geq 0$
for all $j=1, \dots, N$ and for almost all $e \in S^2$, we
suppose, to the contrary, that there is a set $U$ of positive
measure in $S^2$, such that for some vertex $p_j$ in 
$\Gamma_1 \cap \Gamma_2$, whenever $e$ is in $U$, the inequality
$ab<0$
is satisfied, where $a={\rm nlm}_{\Gamma_1}(e, p_j)$ and 
$b={\rm nlm}_{\Gamma_2}(e, p_j)$. Then for $e$ in $U$,
$a$ and $b$ are of opposite signs. Let $U_1$ be the
part of $U$ where $a< 0<b$
holds: we may assume $U_1$ has positive measure, otherwise
exchange $\Gamma_1$ with $\Gamma_2$.  On $U_1$, we have
\[
 (a+b)^+ < b^+= a^+ + b^+.
\]
Recall that $a+b={\rm nlm}_{\Gamma}(e,p_j).$ 
Hence the inequality between half-integers
$${\rm nlm}_\Gamma^+(e,p_j)< 
{\rm nlm}_{\Gamma_1}^+(e,p_j)+{\rm nlm}_{\Gamma_2}^+(e, p_j)$$ 
is valid on the set of positive measure $U_1$, which in
turn implies that
$\nc(\Gamma_1 \cup \Gamma_2) < \nc(\Gamma_1) + \nc(\Gamma_2)$, 
contradicting the assumption of equality.  
\qed

%
\subsection{One-point union of graphs}

%
\begin{pro}\label{1ptunion}
If the graph $\Gamma$ is the one-point union of graphs
$\Gamma_1$ and $\Gamma_2$, where the points $p_1$ chosen in $\Gamma_1$
and $p_2$ chosen in $\Gamma_2$ are not topological vertices,
then the minimum $\nc$ among all mappings is subadditive, and the
minimum $\nc$ minus $2\pi$ is superadditive:
$$\nc(\{\Gamma_1\})+\nc(\{\Gamma_2\})-2\pi\leq \nc(\{\Gamma\})\leq
\nc(\{\Gamma_1\})+\nc(\{\Gamma_2\}).$$
Further, if the points $p_1 \in\Gamma_1$ and $p_2 \in\Gamma_2$ 
may appear as
extreme points on mappings of minimum $\nc$, then the minimum
net total curvature among all mappings, minus $2\pi$, is
additive:
$$\nc(\{\Gamma\})=
\nc(\{\Gamma_1\})+\nc(\{\Gamma_2\})-2\pi.$$
\end{pro}
\pf
Write $p \in \Gamma$ for the identified points $p_1=p_2=p$.

Choose flat 
mappings $f_1:\Gamma_1\to\Re$ and $f_2:\Gamma_2\to\Re$, 
adding constants so that the chosen points 
$p_1 \in \Gamma_1$ and $p_2 \in \Gamma_2$ have
$f_1(p_1)=f_2(p_2)=0.$ Further, by Proposition \ref{minmonot}, we
may assume that $f_1$ and $f_2$ are strictly monotone on the edges
of
$\Gamma_1$ resp. $\Gamma_2$ containing $p_1$ resp. $p_2$.
Let $f:\Gamma\to\Re$ be defined as $f_1$ on $\Gamma_1$
and as $f_2$ on $\Gamma_2$. Then at the common point 
of $\Gamma_1$ and $\Gamma_2$, 
$f(p)=0$, and $f$ is continuous. But since $f_1$ and $f_2$  	
are monotone on the edges containing $p_1$ and $p_2$,		
${\rm nlm}_{\Gamma_1}(p_1)=0={\rm nlm}_{\Gamma_2}(p_2)$,
so we have $\nc(\{\Gamma\})\leq \nc(f)=
\nc(f_1)+ \nc(f_2)= 			
\nc(\{\Gamma_1\})+\nc(\{\Gamma_2\})$
by Proposition \ref{subadditivity}.

Next, for all $g:\Gamma \to \Re$, we shall show 
that ${\rm NTC}(g) \geq
{\rm NTC}(\{\Gamma_1\})+{\rm NTC}(\{\Gamma_2\})-2\pi.$
Given $g$, write $g_1$ resp. $g_2$ for the restriction of $g$ to
$\Gamma_1$ resp. $\Gamma_2$. Then 
$\mu_g(e)=\mu_{g_1}(e)-{\rm nlm}^+_{g_1}(p_1)
+\mu_{g_2}(e) -{\rm nlm}^+_{g_2}(p_2) +{\rm nlm}^+_g(p).$
Now for any real numbers $a$ and $b$, the difference 
$(a+b)^+ - (a^++b^+)$ is equal to $\pm a$, $\pm b$ or $0$,
depending on the various signs. Let $a={\rm nlm}_{g_1}(p_1)$ and 
$b={\rm nlm}_{g_2}(p_2)$. Then 
since $p_1$ and $p_2$ are not topological 	
vertices of $\Gamma_1$ resp. $\Gamma_2$, 	
$a,b \in \{-1,0,+1\}$ and
$a+b = {\rm nlm}_{g}(p)$ by Lemma \ref{combin}. In any case, we
have 
$${\rm nlm}_{g}^+(p)-{\rm nlm}_{g_1}^+(p_1)-
{\rm nlm}_{g_2}^+(p_2) \geq -1.$$
Thus, $\mu_g(e)\geq \mu_{g_1}(e)+\mu_{g_2}(e) -1$, and
multiplying by $2\pi$, 
${\rm NTC}(g) \geq
\nc(g_1)+\nc(g_2)-2\pi \geq
{\rm NTC}(\{\Gamma_1\})+{\rm NTC}(\{\Gamma_2\})-2\pi.$

Finally, assume $p_1$ and $p_2$ are extreme
points for flat mappings $f_1:\Gamma_1 \to \Re$ resp. 
$f_2:\Gamma_2 \to \Re$. 
We may assume that 	
$f_1(p_1)=0=\min f_1(\Gamma_1)$ and 
$f_2(p_2)=0=\max f_2(\Gamma_2)$. Then 
${\rm nlm}_{f_2}(p_2)=1$ and ${\rm nlm}_{f_1}(p_1)=-1$, and hence
using Lemma \ref{combin}, ${\rm nlm}_{f}(p)=0$. So 
$\mu_f(e)=\mu_{f_1}(e)-{\rm nlm}^+_{f_1}(p_1)
+\mu_{f_2}(e)-{\rm nlm^+}_{f_2}(p_2)+{\rm nlm}^+_f(p)=
\mu_{f_1}(e)+\mu_{f_2}(e)-1.$ 
Multiplying by $2\pi$, we have 
${\rm NTC}(\{\Gamma\})\leq {\rm NTC}(f)=
{\rm NTC}(\{\Gamma_1\})+{\rm NTC}(\{\Gamma_2\})-2\pi.$
\qed\\

%
\section{NET TOTAL CURVATURE FOR DEGREE $3$}\label{deg3}
%
%
\subsection{Simple description of net total curvature}\label{simple}

%
\begin{pro}\label{net3}
For any graph $\Gamma$ and any parameterization $\,\Gamma'$ of its
double, $\nc(\Gamma) \leq \frac12 {\mathcal C}(\Gamma')$.
If $\, \Gamma$ is a {\em trivalent} graph, that is, having vertices 
of degree at most three, then 
$\nc(\Gamma) = \frac12 {\mathcal C}(\Gamma')$ 
for any parameterization $\Gamma'$ 
which does not immediately repeat any edge of $\Gamma$.
\end{pro} 
\pf
The first conclusion follows from Corollary \ref{mucompare}.

Now consider a trivalent graph $\Gamma$.
Observe that $\Gamma'$ would be forced to immediately repeat
any edge which ends in a vertex of degree $1$; thus, we may assume
that $\, \Gamma$ has only vertices of degree $2$ or $3$.
Since $\Gamma'$ covers each edge of $\Gamma$ twice, we need only
show, for every vertex $q$ of $\Gamma$, having degree 
$d = d(q) \in \{2,3\}$, that
%
%
\begin{equation}\label{*}
2\, {\rm ntc}_\Gamma(q)= \sum_{i=1}^d {\rm c}_{\Gamma'}(q_i),
\end{equation}
where $q_1, \dots, q_d$ are the vertices of $\Gamma'$ over $q$.
If $d=2$, since $\Gamma'$ does not immediately repeat any edge of
$\Gamma,$ we have
${\rm ntc}_\Gamma(q)={\rm c}_{\Gamma'}(q_1)=
{\rm c}_{\Gamma'}(q_2)$, 
so equation \eqref{*} clearly holds.  For
$d=3$, write both sides of equation \eqref{*} as integrals over
$S^2$, using the definition \eqref{defnc} of ${\rm
ntc}_\Gamma(q)$.
Since $\Gamma'$ does not immediately repeat any edge, the three
pairs of tangent vectors
$\{T_1^{\Gamma'}(q_j),T_2^{\Gamma'}(q_j)\}$, $1\leq j\leq 3$,  
comprise all three pairs taken from the triple
$\{T_1^\Gamma(q),T_2^\Gamma(q),T_3^\Gamma(q)\}$.  We need to show
that
\begin{eqnarray*}
2\int_{S^2}\left[\chi_1+\chi_2+\chi_3\right]^+\,dA_{S^2}&=&
\int_{S^2}\left[\chi_1+\chi_2\right]^+\,dA_{S^2}+ \\
+ \int_{S^2}\left[\chi_2+\chi_3\right]^+\,dA_{S^2}&+&
\int_{S^2}\left[\chi_3+\chi_1\right]^+\,dA_{S^2},
\end{eqnarray*}
where at each direction $e\in S^2$, $\chi_j(e) = \pm 1$ is the
sign of $\langle -e, T_j^\Gamma(q)\rangle$.  But the integrands
are equal at almost every point $e$ of $S^2$:
$$
2\left[\chi_1+\chi_2+\chi_3\right]^+ =
\left[\chi_1+\chi_2\right]^+ +
\left[\chi_2+\chi_3\right]^+ +
\left[\chi_3+\chi_1\right]^+,
$$
as may be confirmed by cases:  $6=6$ if $\chi_1=\chi_2=\chi_3=+1$;
$2=2$ if exactly one of the $\chi_i$ equals $-1$, and $0=0$ in the
remaining cases.
\qed
\vspace{1em}

%
\subsection{Simple description of net total
curvature fails, $d \geq 4$}

%
\begin{obs}\label{notinf}
We have seen in Corollary \ref{mucompare} that for graphs with
vertices of degree $\leq 3$, if a parameterization $\Gamma'$ of
the double $\widetilde\Gamma$ of $\Gamma$ does not immediately
repeat any edge of $\Gamma$, then 
$\nc(\Gamma) = \frac12 {\mathcal C}(\Gamma')$, the total curvature
in the usual sense of the link $\Gamma'$.  
A natural suggestion would be that for
general graphs $\Gamma$, $\nc(\Gamma)$ might be half the infimum
of total curvature of all such parameterizations $\Gamma'$ of the
double.  However, in some cases, we have the {\bf strict
inequality} $\nc(\Gamma) < \inf_{\Gamma'}\frac12 \nc(\Gamma')$.
\end{obs}

In light of Proposition \ref{net3}, we choose an example of a
vertex $q$ of degree four, and consider the local contributions
to $\nc$ for $\Gamma = K_{1,4}$ and for $\Gamma'$, which is
the union of four arcs.

Suppose that for a small positive angle $\alpha$, ($\alpha \leq 1$
radian would suffice) the four unit tangent vectors at $q$ are
$T_1 = (1,0,0)$; $T_2 = (0,1,0)$;
$T_3=(-\cos\alpha,0,\sin\alpha)$;
and $T_4 = (0,-\cos\alpha,-\sin\alpha)$. Write the exterior angles
as $\theta_{ij} = \pi - \arccos \langle T_i, T_j \rangle.$ Then
$\inf_{\Gamma'}\frac12 {\mathcal C}(\Gamma')=
\theta_{13}+\theta_{24} = 2\alpha.$ However, ${\rm ntc}(q)$ is
strictly less than $2\alpha$. This may be seen by writing
${\rm ntc}(q)$ as an integral over $S^2$, according to the
definition \eqref{defnc}, and noting that cancellation occurs
between two of the four lune-shaped sectors.
\qed

%
\subsection{Minimum NTC for trivalent graphs}

\vspace{1em}

Using the relation $\nc(\Gamma) = \frac12 \nc(\Gamma')$ between
the net total curvature of a given trivalent graph $\Gamma$ and
the total curvature for a non-reversing double cover $\Gamma'$ of
the graph, we can determine the minimum net total curvature of a
trivalent graph embedded in $\Re^n$, whose value is then related
to the Euler characteristic of the graph $\chi(\Gamma)=- k/2$.

First we introduce the following definition.

%
%
\begin{defi}\label{bridge}
For a given graph $\Gamma$ and a mapping $f:\Gamma \to \Re$, let
the {\em extended bridge number} $B(f)$ be one-half the number
of local extrema.  Write $B(\{\Gamma\})$ for the minimum of
$B(f)$ among all mappings $f:\Gamma \to \Re$.  For a given isotopy
type $[\Gamma]$ of embeddings into $\Re^3$, let $B([\Gamma])$ be
one-half the minimum number of local extrema for a mapping
$f:\Gamma\to\Re$ in the closure of the isotopy class $[\Gamma]$.  

\end{defi}

For an integer $m\geq 3$, let $\theta_m$ be the graph with two
vertices $q^+, q_-$ and $m$ edges, each of which has $q^+$ and
$q^-$ as its two endpoints. Then $\theta = \theta_3$ has the form
of the lower-case Greek letter $\theta$.  

%
\begin{rem}
For a knot, the number of local maxima equals the number of local
minima. The minimum number of local maxima is called the 
{\em bridge number}, and equals the number of local minima. This is 
consistent with our Definition \ref{bridge} of the extended bridge 
number.  Of course, for knots, the minimum bridge number among all
isotopy classes $B(\{S^1\})=1$, and only $B([S^1])$ is of
interest for a specific isotopy class $[S^1]$. For certain graphs,
the minimum numbers of local maxima and local minima may not occur
at the same time for any mapping: see the example of Observation
\ref{B>1} below. For isotopy classes of $\theta$-graphs, Goda
\cite{Go} has given a definition of an integer-valued bridge index
which is similar in spirit to the definition above.  
\end{rem}

%
\begin{thm}\label{trivalent}
If $\, \Gamma$ is a trivalent graph, and if $f_0:\Gamma\to\Re$ is
monotone on topological edges and has the minimum number
$2B(\{\Gamma\})$ of local extrema, then
$\nc(f_0) = \nc(\{\Gamma\})= 
\pi\Big(2B(\{\Gamma\}) + \frac{k}{2}\Big)$,
where $k$ is the number of topological vertices of $\, \Gamma.$
For a given isotopy class $[\Gamma]$, 
$\nc([\Gamma])= \pi\Big(2B([\Gamma]) + \frac{k}{2}\Big)$.
\end{thm}

\pf
Recall that $\nc(\{\Gamma\})$ denotes the infimum of $\nc(f)$
among $f:\Gamma \to \Re^3$ or among $f:\Gamma \to \Re$, as may be
seen from Corollary \ref{1dsuff}.  

We first consider a mapping $f_1:\Gamma \to \Re$ with the property
that any local maximum or local minimum points of $f_1$ are
interior points of topological edges. Then all topological
vertices $v$, since they have degree $d(v)=3$ and $d^\pm(v)\neq
0$, 
have ${\rm nlm}(v)=\pm 1/2$, by Proposition \ref{val3}.
Let $\Lambda$ be the number
of local maximum points of $f_1$, $V$ the number of local minimum
points, $\lambda$ the number of vertices with ${\rm nlm}=+1/2$,
and ${\tt y}$ the number of vertices with ${\rm nlm}=-1/2$.
Then $\lambda + {\tt y} = k$, the total number of
vertices, and $\Lambda + V \geq 2B(\{\Gamma\})$. Hence applying 
Corollary \ref{absnlm},
%
\begin{equation}\label{exacttri}
\mu = \frac12\sum_v|{\rm nlm}(v)|=
\frac12[\Lambda+V+\frac{\lambda+{\tt y}}{2}]\geq B(\{\Gamma\})+k/4,
\end{equation}
with equality iff $\Lambda + V = 2B(\{\Gamma\})$.

We next consider any mapping $f_0:\Gamma \to \Re$ in general
position: in particular, the critical valuess of $f_0$ are
isolated. In a similar fashion to the proof of Proposition
\ref{minmonot}, we shall replace $f_0$ with a mapping whose local
extrema are not topological vertices. Specifically, if $f_0$
assumes a local maximum at any topological vertex $v$, then, since
$d(v)=3$, ${\rm nlm}_{f_0}(v)=3/2$. $f_0$ may be isotoped in a
small neighborhood of $v$ to $f_1:\Gamma \to \Re$ so that near
$v$, the local maximum occurs at an interior point $q$ of one of
the three edges with endpoint $v$, and thus 
${\rm nlm}_{f_1}(q)=1$; while the up-degree $d_{f_1}^+(v)=1$ and
the down-degree $d_{f_1}^-(v)=2$, so that ${\rm nlm}_{f_1}(v)$ is 
now $\frac12$.  Thus, $\mu_{f_1}(e)=\mu_{f_0}(e)$.  Similarly, if
$f_0$ assumes a local minimum at a topological vertex $w$, then
$f_0$ may be isotoped in a neighborhood of $w$ to 
$f_1:\Gamma \to \Re$ so that the local minimum of $f_1$ near $w$
occurs at an interior point of any of the three edges with
endpoint $w$, and $\mu_{f_1}(e)=\mu_{f_0}(e)$.  Then 
any local
extreme points of $f_1$ are interior points of topological edges.
Thus, we have shown that $\mu_{f_0}(e) \geq B(\{\Gamma\}) + k/4$,
with equality if $f_1$ has exactly $2B(\{\Gamma\})$ as its number
of local extrema, which holds iff $f_0$ has the minimum number
$2B(\{\Gamma\})$ of local extrema.

Thus 
$\nc(\{\Gamma\})=2\pi\mu_{f_0}(e)=
2\pi\Big(B(\{\Gamma\})+k/4\Big)=
\pi\Big(2B(\{\Gamma\})+k/2\Big).$

Similarly, for a given isotopy class $[\Gamma]$ of embeddings into
$\Re^3$, we may choose $f_0:\Gamma\to\Re$ in the closure of the
isotopy class, deform $f_0$ to a mapping $f_1$ in the closure of
$[\Gamma]$ having no topological vertices as local extrema and
count $\mu_{f_0}(e)=\mu_{f_1}(e)\geq B([\Gamma]) + k/4$, with
equality if $f_0$ has the minimum number $2B([\Gamma])$ of local
extrema. This shows that
$\nc([\Gamma])=\pi\Big(2B([\Gamma])+k/2\Big).$
\qed

%
\begin{rem}
An example geometrically illustrating the lower bound is given by
the dual graph $\Gamma^*$ of the one-skeleton $\Gamma$ of a
triangulation of $S^2$, with the $\{\infty\}$ not coinciding with
any of the vertices of $\Gamma^*$.  The Koebe-Andreev-Thurston
theorem says that there is a circle packing which realizes the
vertex set of $\Gamma^*$ as the set of centers of the circles 
(see \cite{S}).			
The so realized $\Gamma^*$, stereographically projected to
$\Re^2 \subset \Re^3$, attains the lower bound of Theorem
\ref{trivalent} with $B(\{\Gamma^*\}) = 1$, namely 
$\nc([\Gamma])=\pi(2+\frac{k}{2})=\pi( 2-\chi(\Gamma^*))$, where
$k$ is the number of vertices.
\end{rem}

%
\begin{cor}\label{muformula}
If $\Gamma$ is a trivalent graph with $k$ topological vertices,
and $f_0:\Gamma \to \Re$ is a mapping in general position, having
$\Lambda$ local maximum points and $V$ local minimum points, then
$$ \mu_{f_0}(e)=\frac12(\Lambda + V) +\frac{k}{4}\geq
B(\{\Gamma\})+\frac{k}{4}.$$ 
\end{cor}
\pf
Follows immediately from the proof of Theorem \ref{trivalent}:
$f_0$ and $f_1$ have the same number of local maximum or minimum
points.
\qed \\

An interesting trivalent
graph is $L_m$, the ``ladder of $m$ rungs" obtained from two unit
circles in parallel planes by adding $m$ line segments (``rungs")
perpendicular to the planes, each joining one vertex on the first
circle to another vertex on the second circle. 
For example, $L_4$ is the $1$-skeleton of the cube in $\Re^3$.
Note that $L_m$ may be embedded in
$\Re^2$, 
and that the bridge number $B(\{L_m\})=1$.  Since $L_m$ has $2m$
trivalent vertices, we may apply Theorem \ref{trivalent} to
compute the minimum $\nc$ for the type of $L_m$:

%
\begin{cor} 
The minimum net total curvature 
$\nc(\{L_m\})$ for graphs of the type of $L_m$ equals $\pi(2+m)$.
\end{cor}

%
\begin{obs}\label{B>1}
For certain connected trivalent graphs $\Gamma$ containing cut points,
the minimum extended bridge number $B(\{\Gamma\})$ may be greater than $1$.
\end{obs}

{\it Example:}
Let $\Gamma$ be the union of three disjoint circles
$C_1, C_2, C_3$ with three edges $E_i$ connecting a point 
$p_i \in C_i$ with a fourth vertex $p_0$, which is not in any of
the $C_i$, and which is a {\it cut point} of $\Gamma$: 
the number of connected components of
$\Gamma\backslash p_0$ is greater than for $\Gamma$. Given 
$f:\Gamma \to \Re$, after a permutation of
$\{1,2,3\}$, we may assume there is a minimum point 
$q_1\in C_1\cup E_1$ and a maximum point $q_3\in C_3\cup E_3$.
If $q_1$ and $q_3$ are both in $C_1\cup E_1$, we may choose 
$C_2$ arbitrarily in what follows.
Restricted to the closed set $C_2 \cup E_2$, $f$ assumes either a
maximum or a minimum at a point $q_2 \neq p_0$. Since 
$q_2 \neq p_0$, $q_2$ is also a local maximum or a local minimum
for 
$f$ on $\Gamma$. That is, $q_1, q_2, q_3$ are all local extrema.
In
the notation of the proof of Theorem \ref{trivalent}, we have the
number of local extrema $V + \Lambda \geq 3$. Therefore
$B(\{\Gamma\}) \geq \frac{3}{2}$, and 
${\rm NTC}(\{\Gamma\})\geq \pi(3+k/2)=5\pi.$

The reader will be able to construct similar trivalent examples
with $B(\{\Gamma\})$ arbitrarily large.
\qed

In contrast to the results of Theorem \ref{trivalent} and of
Theorem \ref{allbut1tri}, below, for trivalent or nearly trivalent
graphs, the minimum of $\nc$ for a given graph type cannot be
computed merely by counting vertices, but depends in a more subtle
way on the topology of the graph:

%
\begin{obs}\label{samedegree}
When $\Gamma$ is not trivalent, the minimum $\nc(\{\Gamma\})$ of
net total curvature for a connected graph $\Gamma$ with
$B(\{\Gamma\})=1$ is not determined by the number of vertices
and their degrees.
\end{obs}

{\it Example:} We shall construct two planar graphs $S_m$ and $R_m$
having the same number of vertices, all of degree $4$.

Choose an integer $m\geq 3$ and take the image of the embedding
$f_\varepsilon$ of the ``sine wave" $S_m$ to be
the union of the polar-coordinate graphs $C_\pm\subset\Re^2$ of two
functions: $r=1 \pm \varepsilon\sin(m\theta)$. 
$S_m$ has $4m$ edges; 
and $2m$ vertices, all of degree $4$, at $r=1$ and 
$\theta = \pi/m, 2\pi/m, \dots , 2\pi$. For 
$0<\varepsilon <1$,
$f_\varepsilon(S_m)=C_+\cup C_-$ is the union of two smooth cycles.  
For small positive $\varepsilon$, $C_+$ and $C_-$ are convex.
The $2m$ vertices all have ${\rm nlm}(q)=0$, so
$\nc(f_\varepsilon)=\nc(C_+)+\nc(C_-)=2\pi+2\pi$.
Therefore $\nc(\{S_m\}) \leq \nc(f_\varepsilon) = 4\pi$.

For the other graph type, let the ``ring graph" $R_m\subset\Re^2$
be constructed by adding $m$ disjoint small circles $C_i$, each crossing
one large circle $C$ 
at two points $v_{2i-1}, v_{2i}$, $1\leq i\leq m$.
Then $R_m$ has $4m$ edges.  We construct $R_m$ so that the $2m$
vertices $v_1,v_2,\dots,v_{2m}$, 
appear in cyclic order around $C$.
Then $R_m$ has the same
number $2m$ of vertices as does $S_m$, all of degree $4$. 
At each vertex $v_j$, we have ${\rm nlm}(v_j)=0$, so in this
embedding, $\nc(R_m) = 2\pi(m+1)$. We shall show that 
$\nc(f_1) \geq 2\pi m$ for any $f_1:R_m \to \Re^3$. According to 
Corollary \ref{1dsuff}, it is enough to show for every 
$f:R_m \to \Re$ that $\mu_{f}\geq m$. We may assume $f$ is monotone
on each topological edge, according to Proposition
\ref{minmonot}. Depending on the order of $f(v_{2i-2}),
f(v_{2i-1})$ and
$f(v_{2i})$, ${\rm nlm}(v_{2i-1})$ might equal $\pm 1$ or
$\pm 2$, but cannot be $0$, as follows from Lemma \ref{combin},
since the unordered pair $\{d^-(v_{2i-1}),d^+(v_{2i-1})\}$ 
may only be $\{1,3\}$ or $\{0,4\}$. Similarly, $v_{2i}$ is
connected by three edges to
$v_{2i-1}$ and by one edge to $v_{2i+1}$. For the same reasons,
${\rm nlm}(v_{2i})$ might equal $\pm 1$ or $\pm 2$, and cannot
$=0$. So $|{\rm nlm}(v_j)| \geq 1$, $1 \leq j \leq 2m$, and thus
by Corollary \ref{absnlm}, $\mu = \frac12\sum_j |{\rm
nlm}(v_j)|\geq m$. Therefore the minimum of net total curvature
$\nc(\{R_m\})\geq 2m\pi$, which is greater than
$\nc(\{S_m\})\leq 4\pi$, since $m\geq 3$.

(A more detailed analysis shows that $\nc(\{S_m\})= 4\pi$ and
$\nc(\{R_m\})= 2\pi(m+1)$.)			
\qed

\vspace{1em}

Finally, we may extend the methods of proof for Theorem
\ref{trivalent} to allow {\bf one} vertex of higher degree:

%
\begin{thm}\label{allbut1tri}
If $\, \Gamma$ is a graph with one vertex $w$ of degree
$d(w)=m \geq 3$, all other vertices being trivalent, 
and if $w$ shares edges with $m$ distinct trivalent vertices, 
then
$\nc(\{\Gamma\})= \pi\Big(2B(\{\Gamma\}) + \frac{k}{2}\Big)$,
where $k$ is the number of vertices of $\, \Gamma$ having odd
degree. For a given isotopy class $[\Gamma]$,
$\nc([\Gamma])\geq \pi\Big(2B([\Gamma]) + \frac{k}{2}\Big)$.
\end{thm}
\pf
Consider any mapping $g:\Gamma \to \Re$ in general position. If
$m$ is even, then $|{\rm nlm}_g(w)|\geq 0$; if $m$ is odd, then
$|{\rm nlm}_g(w)|\geq \frac{1}{2}$, by Proposition \ref{val3}. 
If some topological vertex is a local
extreme point, then as in the proof of Theorem \ref{trivalent},
$g$ may be modified without changing $\nc(g)$ so that all
$\Lambda+V$ local extreme points are interior points of edges,
with ${\rm nlm}=\pm 1$.  By Corollary \ref{absnlm}, we have 
$\mu_g(e)=\frac12\sum|{\rm nlm}(v)|\geq
\frac12\Big(\Lambda+V+\frac{k}{2}\Big)\geq
B(\{\Gamma\})+\frac{k}{4}$. This shows that
$$\nc(\{\Gamma\})\geq \pi\Big(2B(\{\Gamma\})+\frac{k}{2}\Big).$$

Now let $f_0:\Gamma \to \Re$ be monotone on topological edges and
have the minimum number $2B(\{\Gamma\})$ of local extreme points
(see Corollary \ref{minmonot}).		
As in the proof of Theorem \ref{trivalent}, $f_0$ may be modified
without changing $\nc(f_0)$ so that all $2B(\{\Gamma\})$ local
extreme points are interior points of edges.  $f_0$ may be further
modified so that the distinct vertices $v_1, \dots, v_m$ which share 
edges with $w$ are balanced: $f(v_j)<f(w)$ for half of the
$j=1,\dots,m$, if $m$ is even, or for half of $m+1$, if $m$ is odd.  
Having chosen $f(v_j)$, we define $f$ along the (unique) edge from
$w$ to $v_j$ to be monotone, for $j=1,\dots,m$.  Therefore if $m$
is even, then ${\rm nlm}_f(w)= 0$; and if $m$ is odd, then 
${\rm nlm}_f(w)= \frac{1}{2}$, by Lemma \ref{combin}.  We compute
$\mu_f(e)=\frac12\sum|{\rm nlm}(v)|=
\frac12(\Lambda+V+\frac{k}{2})= B(\{\Gamma\})+\frac{k}{4}$.  We
conclude that 
$\nc(\{\Gamma\})= \pi\Big(2B(\{\Gamma\})+\frac{k}{2}\Big)$.

For a given isotopy class $[\Gamma]$, the proof is analogous to
the above. Choose a mapping $g:\Gamma \to \Re$ in the closure of
$[\Gamma]$, and modify $g$ without leaving the closure of
the isotopy class. Choose $f:\Gamma \to \Re$ which has the minimum
number $2B([\Gamma])$ of local extreme points, and modify it so
that topological vertices are not local extreme points. In contrast 
to the proof of Theorem \ref{trivalent}, a balanced
arrangement of vertices may not be possible in the given isotopy
class.  In any case, if $m$ is even, then 
$|{\rm nlm}_f(w)|\geq 0$; and if $m$ is odd, 
$|{\rm nlm}_f(w)|\geq \frac{1}{2}$, by Proposition \ref{val3}.
Thus applying Corollary \ref{absnlm}, we find 
$\nc([\Gamma])\geq\pi\Big(2B([\Gamma])+\frac{k}{2}\Big)$.
\qed

%
\begin{obs}
When all vertices of $\Gamma$ are trivalent except $w$, $d(w)\geq 4$,
and when $w$ shares more than one edge with another vertex of $\Gamma$,
then in certain cases, 
$\nc(\{\Gamma\})>\pi\Big(2B(\{\Gamma\})+\frac{k}{2}\Big)$, where $k$
is the number of vertices of odd degree.
\end{obs}
{\it Example:} Choose $\Gamma$ to be the one-point union of
$\Gamma_1$, $\Gamma_2$ and $\Gamma_3,$ where $\Gamma_i=
\theta=\theta_3$,
$i=1,2,3$, and the point $w_i$ chosen from $\Gamma_i$ is one of
its two vertices $v_i, w_i$. Then the identified point
$w=w_1=w_2=w_3$ of $\Gamma$ has $d(w)=9$, and each of the other
three vertices $v_1, v_2, v_3$ has degree $3$.  

Choose a flat
map $f:\Gamma\to\Re$. We may assume that $f$ is monotone on each
edge, applying Proposition \ref{minmonot}. If
$f(v_1)<f(v_2)<f(w)<f(v_3)$, then $d^+(w)=3$, $d^-(w)=6$, so 
${\rm nlm}(w)=\frac32$, while $v_i$ is a local extreme point, so
${\rm nlm}(v_i)=\pm\frac32$, $1=1,2,3$.  This gives $\mu = 3$.
The case where $f(v_1)<f(w)<f(v_2)<f(v_3)$ is similar.  If $w$ is
an extreme point of $f$, then ${\rm nlm}(w)=\pm\frac92$ and
$\mu\geq\frac92>3$, contradicting flatness of $f$. This shows that 
$\nc(\{\Gamma\})=\nc(f)=6\pi$.

On the other hand, we may show as in Observation \ref{B>1} that
$B(\{\Gamma\})=\frac32$. All four vertices have odd degree, so
$k=4$, and $\pi\Big(2B(\{\Gamma\})+\frac{k}{2}\Big) = 5\pi$.
\qed

\vspace{1em}

Let $W_m$ denote the ``wheel" of $m$ spokes, consisting of a cycle
$C$ containing $m$ vertices $v_1,\dots,v_m$ (the ``rim"), a
central vertex $w$ (the ``hub") not on $C$, and edges $E_i$ (the
``spokes") connecting $w$ to $v_i$, $1\leq i \leq m$.

%
\begin{cor}
The minimum net total curvature $\nc(\{W_m\})$ for graphs in
$\Re^3$ homeomorphic to $W_m$ equals
$\pi(2+\lceil\frac{m}{2}\rceil)$.
\end{cor}
\pf
We have one ``hub" vertex $w$ with $d(w)=m$, and all other  	
vertices have degree $3$. Observe that the bridge number
$B(\{W_m\})=1$.  According to Theorem \ref{allbut1tri}, we have
$\nc(\{W_m\})= \pi\Big(2B(\{W_m\}) + \frac{k}{2}\Big)$, where $k$
is the number of vertices of odd degree: $k=m$ if $m$ is even, or
$k=m+1$ if $m$ is odd: $k=2\lceil\frac{m}{2}\rceil$. Thus 
$\nc(\{W_m\})= \pi\Big(2+\lceil\frac{m}{2}\rceil\Big)$.
\qed\\

%
\section{LOWER BOUNDS OF NET TOTAL CURVATURE}\label{lowbds}

The {\em width} of an isotopy class $[\Gamma]$ of embeddings 	
of a graph $\Gamma$ into $\Re^3$ is the minimum among 		
representatives of the class of the maximum number of points	
of the graph meeting a family of parallel planes. 		
More precisely, we write
${\rm width}([\Gamma]):=
\min_{f:\Gamma \to \Re^3\vert f\in[\Gamma]} 
\min_{e\in S^2} \max_{s\in\Re} \#(e,s).$ 
For any homeomorphism type $\{\Gamma\}$ define 
${\rm width}(\{\Gamma\})$ to be the minimum over isotopy types.

\begin{thm}\label{incrdecr}
Let $\Gamma$ be a graph, and consider an isotopy class $[\Gamma]$
of embeddings $f:\Gamma \to \Re^3$. Then 
$$\nc([\Gamma]) \geq \pi \ {\rm width}([\Gamma]).$$ 
As a consequence, 
$\nc(\{\Gamma\}) \geq \pi \ {\rm width}(\{\Gamma\}).$
Moreover, if for some $e\in S^2$, an embedding 
$f:\Gamma \to \Re^3$ and $s_0\in \Re$, the integers
$\#(e,s)$ are increasing in $s$ for $s<s_0$ and decreasing for
$s>s_0$, then $\nc([\Gamma]) = \#(e,s_0)\,\pi.$
\end{thm}
\pf
Choose an embedding $g:\Gamma \to \Re^3$ in the given isotopy
class, with $\max_{s\in\Re} \#(e,s)={\rm width}([\Gamma])$. 
There exist $e\in S^2$ and $s_0\in\Re$ with 
$\#(e,s_0) = \max_{s\in\Re} \#(e,s) = {\rm width}([\Gamma])$.
Replace $e$ if necessary by a nearby point in
$S^2$ so that the values $g(v_i)$, $i=1,\dots,m$ are distinct.
Next do cylindrical shrinking: without changing
$\#(e,s)$ for $s\in \Re$, shrink the image of $g$ in directions
orthogonal to $e$ by a factor $\delta>0$ to obtain a family 
$\{g_\delta\}$ from the same isotopy class $[\Gamma]$, with
$\nc(g_\delta)\to \nc(g_0)$, where we may identify 
$g_0: \Gamma \to \Re e \subset\Re^3$ with 
$p_e\circ g=p_e\circ g_\delta:\Gamma \to \Re$.  But
$$\nc(p_e\circ g)=
\frac12\int_{S^2}\mu(u)\,dA_{S^2}(u)=2\,\pi\,\mu(e),$$
since for $p_u\circ p_e\circ g$, the local maximum and minimum
points are the same as for $p_e\circ g$ if $\langle e,u\rangle >0$
and reversed if $\langle e,u\rangle <0$ (recall that
$\mu(-e)=\mu(e)$).

We write the topological vertices and the local extrema of $g_0$
as $v_1,\dots,v_m$.  
Let the indexing be chosen so that $g_0(v_i)<g_0(v_{i+1})$,
$i=1,\dots, m-1$.  Now estimate $\mu(e)$ from below: using Lemma
\ref{combin},
%
\begin{equation}\label{mueqno} 
\mu(e)= 
\sum_{i=1}^m {\rm nlm}^+_{g_0}(v_i)\geq 
\sum_{i=k+1}^m {\rm nlm}_{g_0}(v_i) 
= \frac12 \#(e,s) 
\end{equation}
for any $s$, $g_0(v_k)<s<g_0(v_{k+1})$.  This shows that
$\mu(e)\geq \frac12 {\rm width}([\Gamma])$, and therefore 
$$\nc(g)\geq\nc(g_0)=2\pi\,\mu(e)\geq\,\pi\,{\rm width}([\Gamma]).$$

Now suppose that the integers $\#(e,s)$ are increasing in $s$ for
$s<s_0$ and decreasing for $s>s_0$. Then for $g_0(v_i)>s_0$, we
have ${\rm nlm}(g_0(v_i))\geq 0$ by Lemma \ref{combin}, and the 
inequality \eqref{mueqno} becomes equality at $s=s_0$.
\qed

\vspace{1em}

\begin{lem}\label{widthK_m}
For an integer $\ell$, the minimum width of the complete graph 
$K_{2\ell}$ on $2\ell$ vertices
is ${\rm width}(\{K_{2\ell}\})=\ell^2$; for $2\ell + 1$ vertices,
${\rm width}(\{K_{2\ell+1}\})= \ell(\ell+1).$
\end{lem}
\pf 
Write $E_{ij}$ for the edge of $K_m$ joining $v_i$ to $v_j$,
$1\leq i < j \leq m$, and suppose $g:K_m\to \Re$ has distinct
values at the vertices: $g(v_1)<g(v_2)<\cdots <g(v_m)$.

Then for any $g(v_k) <s< g(v_k+1)$, there are $k(m-k)$ edges
$E_{ij}$ with $i\leq k <j$; each of these edges has at least
one interior point mapping to $s$, which shows that 
$\#(e,s)\geq k(m-k).$ If $m$ is even: $m=2\ell$, these lower bounds
have the maximum value $\ell^2$ when  $k=\ell$. If $m$ is odd:
$m=2\ell+1,$ these lower bounds have the maximum value
$\ell(\ell+1)$ when $k=\ell$ or $k=\ell+1$. This shows that the
width of $K_{2\ell} \geq\ell^2$ and the width of $K_{2\ell+1}\geq
\ell(\ell+1).$ On the other hand, equality holds for the piecewise
linear embedding of $K_m$ into $\Re$ with vertices in general
position and straight edges $E_{ij}$, which shows that 
${\rm width}(\{K_{2\ell}\}) = \ell^2$ and 
${\rm width}(\{K_{2\ell+1}\}) = \ell(\ell+1).$
\qed

\vspace{1em}

%
\begin{pro}\label{NTCK_m}
For all $g:K_m\to \Re$, $\nc(g) \geq \pi\,\ell^2$ if $m=2\ell$ is
even; and $\nc(g) \geq \pi\,\ell(\ell+1)$ if $m=2\ell+1$ is odd.
Equality holds for an embedding of $K_m$ into $\Re$ with vertices
in general position and monotone on each edge; therefore
$\nc(\{K_{2\ell}\})=\pi\,\ell^2$, and
$\nc(\{K_{2\ell+1}\})=\pi\,\ell(\ell+1)$.
\end{pro}
\pf
The lower bound on $\nc(\{K_m\})$ follows from Theorem
\ref{incrdecr} and Lemma \ref{widthK_m}.

Now suppose $g:K_m\to \Re$ is monotone on each edge, and number
the vertices of $K_m$ so that for all $i$, $g(v_i)< g(v_{i+1})$.
Then as in the proof of Lemma \ref{widthK_m}, $\#(e,s)=k(m-k)$ for 
$g(v_k)<s<g(v_{k+1})$. These cardinalities are increasing for 
$0\leq k \leq\ell$ and decreasing for $\ell+1<k<m$. Thus, if 
$g(v_\ell)<s_0<g(v_{\ell+1})$, then by Theorem \ref{incrdecr}, 
$\nc([\Gamma]) = \#(e,s_0)\,\pi= \ell(m-\ell)\,\pi,$ as claimed.
\qed\\

Let $K_{m,n}$ be the complete bipartite graph with $m+n$ vertices
divided into two sets: $v_i, 1\leq i \leq m$ and 
$w_j, 1\leq j\leq n$, having one edge $E_{ij}$ joining $v_i$ to
$w_j$, for each $1\leq i \leq m$ and $1\leq j\leq n$.

%
\begin{pro}\label{bipartite}
$\nc(\{K_{m,n}\})=\lceil\frac{mn}{2}\rceil\,\pi$.
\end{pro}
\pf
$K_{m,n}$ has vertices $v_1,\dots, v_m$ of degree $d(v_i)=n$ and
vertices $w_1,\dots, w_n$ of degree $d(w_j)=m$. Consider a mapping
$g:K_{m,n} \to \Re$ in general position, so that the $m+n$
vertices of $K_{m,n}$ have distinct images. We wish to show
$\mu(e)=\mu_g(e)\geq\frac{mn}{4},$ if $m$ or $n$ is even, or
$\frac{mn+1}{4},$ if both $m$ and $n$ are odd. 

For this purpose, according to Proposition \ref{minmonot},
we may first reduce $\mu(e)$ or leave it unchanged by replacing
$g$ with a mapping (also called $g$) which is monotone on each
edge $E_{ij}$ of $K_{m,n}$. The values of ${\rm nlm}(w_j)$ and
of ${\rm nlm}(v_i)$ are now determined by the order of the vertex
images $g(v_1),\dots,g(v_m),g(w_1),\dots,g(w_n)$. Since $K_{m,n}$
is symmetric under permutations of $\{v_1,\dots,v_m\}$ and
permutations of $\{w_1,\dots,w_n\}$, we shall assume that
$g(v_i)<g(v_{i+1})$, $i=1,\dots,m-1$ and $g(w_j)<g(w_{j+1})$,
$j=1,\dots,n-1$. For $i=1,\dots, m$ we write $k_i$ for the largest
index $j$ such that $g(w_j)<g(v_i).$ Then 
$0\leq k_1\leq \dots\leq k_m \leq n$, and these integers determine 
$\mu(e)$. According to
Lemma \ref{combin}, ${\rm nlm}(v_i)=k_i-\frac{n}{2}, i=1,\dots,m$.
For $j\leq k_1$ and for $j\geq k_m+1$, we have 
${\rm nlm}(w_j)=\pm\frac{m}{2}$; for $k_1<j\leq k_2$ and for
$k_{m-1}<j\leq k_m$, we find 
${\rm nlm}(w_j)=\pm\Big(\frac{m}{2}-1\Big)$; and so on until we
find ${\rm nlm}(w_j)=0$ on the middle interval 
$k_p<j\leq k_{p+1}$, if $m=2p$ is even; or, if $m=2p+1$ is odd,
${\rm nlm}(w_j)=-\frac12$ for $k_p<j\leq k_{p+1}$ and 
${\rm nlm}(w_j)=+\frac12$ for the other middle interval
$k_{p+1}<j\leq k_{p+2}$. Thus according to Lemma \ref{combin} and
Corollary \ref{absnlm}, if $m=2p$ is {\bf even},
%
%
\begin{eqnarray}\label{muformeven} 
2\mu(e)&=&\sum_{i=1}^m|{\rm nlm}(v_i)|+\sum_{j=1}^n|{\rm
nlm}(w_j)|=
\sum_{i=1}^m|k_i-\frac{n}{2}|+(k_1+n-k_m)\frac{m}{2}\nonumber\\ 
&+& (k_2-k_1+k_m-k_{m-1})\Big[\frac{m}{2}-1\Big]+
\dots \nonumber\\
 &+& (k_p-k_{p-1}+k_{p+2}-k_{p+1})
\Big[\frac{m}{2}-(p-1)\Big] + (k_{p+1}-k_p)\Big[0\Big]\\ 
 &=& 
\sum_{i=1}^m\Big|k_i-\frac{n}{2}\Big|+\frac{mn}{2}+\sum_{i=1}^p
k_i-
\sum_{i=p+1}^m k_i \nonumber\\
&=&
\frac{mn}{2}+\sum_{i=1}^p\Big[|k_i-\frac{n}{2}|+(k_i-\frac{n}{2})\Big]+
\sum_{i=p+1}^m\Big[|k_i-\frac{n}{2}|-(k_i-\frac{n}{2})\Big].\nonumber 
\end{eqnarray}
Note that formula \eqref{muformeven} assumes its minimum value
$2\mu(e)=\frac{mn}{2}$ when 
$k_1\leq\dots\leq k_p\leq\frac{n}{2}\leq k_{p+1}\leq\dots\leq
k_m.$

If $m=2p+1$ is {\bf odd}, then
%
%
\begin{eqnarray}\label{muformodd}
2\mu(e)&=& \sum_{i=1}^m|k_i-\frac{n}{2}|+(k_1+n-k_m)\frac{m}{2}+
(k_2-k_1+k_m-k_{m-1})\Big[\frac{m}{2}-1\Big]+
\dots \nonumber\\
&+&(k_{p+3}-k_{p+2})\Big[\frac{m}{2}-(p-1)\Big]+  
(k_{p+2}-k_p)\Big[\frac12\Big]= \nonumber \\
 &=&\sum_{i=1}^m|k_i-\frac{n}{2}|+\frac{mn}{2}+
\sum_{i=1}^p k_i- \sum_{i=p+2}^m k_i\\
&=&
\frac{mn}{2}+\sum_{i=1}^p\Big[|k_i-\frac{n}{2}|+(k_i-\frac{n}{2})\Big]+
\sum_{i=p+2}^m\Big[|k_i-\frac{n}{2}|-(k_i-\frac{n}{2})\Big]+
|k_{p+1}-\frac{n}{2}|.\nonumber
\end{eqnarray}
Observe that formula \eqref{muformodd} has the minimum value
$2\mu(e)=\frac{mn}{2}$ when $n$ is even and 
$k_1\leq\dots\leq k_p\leq\frac{n}{2}= k_{p+1}\leq\dots\leq k_m.$
If $n$ as well as $m$ is odd, then the last term
$|k_{p+1}-\frac{n}{2}|\geq \frac12$, and the minimum value of
$2\mu(e)$ is $\frac{mn+1}{2}$, attained iff 
$k_1\leq\dots\leq k_p\leq\frac{n}{2}\leq k_{p+2}\leq\dots\leq
k_m.$

This shows that for either parity of $m$ or of $n$,
$\mu(e)\geq\frac{mn}{4}$.  If $n$ and $m$ are both odd, we have
the stronger inequality $\mu(e)\geq\frac{mn+1}{4}$.  We may
summarize these conclusions as
$2\mu(e)\geq\lceil\frac{mn}{2}\rceil$, and therefore 
as in the proof of Corollary \ref{1dsuff},
$\nc(\{K_{m,n}\})\geq\lceil\frac{mn}{2}\rceil\,\pi,$ as
we wished to show. 

By abuse of notation, write the formula \eqref{muformeven} or
\eqref{muformodd} as $\mu(k_1,\dots,k_m)$.

To show the inequality in the opposite direction,  	
we need to find a mapping
$f:K_{m,n}\to \Re$ with $\nc(f)=\frac{mn\,\pi}{2}$ ($m$ or $n$ even)
or $\nc(f)=\frac{(mn+1)\,\pi}{2}$ ($m$ and $n$ odd). The above
computation suggests choosing $f$ with $f(v_1), \dots, f(v_m)$
together in the middle of the images of the $w_j$. Write $n=2\ell$
if $n$ is even, or $n=2\ell+1$ if $n$ is odd. Choose values
$f(w_1)<\dots<f(w_\ell)<f(v_1)<\dots<f(v_m)<f(w_{\ell+1})<\dots<f(w_n)$,
and extend $f$ monotonically to each of the $mn$ edges $E_{ij}$.
From formulas \eqref{muformeven} and \eqref{muformodd}, we have
$\mu_f(e)=\mu(\ell,\dots,\ell)=\frac{mn}{4}$, if $m$ or $n$ is
even; or $\mu_f(e)=\mu(\ell,\dots,\ell)=\frac{mn+1}{4}$, if $m$
and $n$ are odd.
\qed

\vspace{1em}

Recall that $\theta_m$ is the graph with two vertices $q^+, q_-$
and $m$ edges.

%
\begin{cor}\label{theta_m}
$\nc(\{\theta_m\})=m\,\pi.$
\end{cor}
\pf
$\theta_m$ is homeomorphic to the complete bipartite graph
$K_{m,2}$, and by the proof of Proposition \ref{bipartite}, we find
$\mu(e)\geq \frac{m}{2}$ for a.a. $e\in S^2$, and hence 
$\nc(\{K_{m,2}\})=m\,\pi$.
\qed

%
%
\section{F\'{A}RY-MILNOR TYPE ISOTOPY
CLASSIFICATION}\label{FaryMilnor}

Recall the F\'{a}ry-Milnor theorem, which states that if the total
curvature of a Jordan curve $\Gamma$ in $\Re^3$ is less than or
equal to $4 \pi$, then $\Gamma$ is unknotted. As we have demonstrated
above, there are a collection of graphs whose values of the
minimum total net curvatures are known.  It is natural to hope
when the net total curvature is small, in the sense of being in a
specific interval to the right of 
the minimal value, that the isotopy type of the graph is restricted,
as is the case for knots: $\Gamma=S^1$. The following 
proposition and corollaries,
however, tell us that results of the F\'{a}ry-Milnor type 
{\bf cannot} be expected to hold for more general graphs.

%
\begin{pro}\label{Gamma_q} 
If $\,\Gamma$ is a graph in $\Re^3$ and if 
$C \subset \Gamma$ is a cycle, 
such that for some $e \in S^2$, $p_e \circ C$ has at least two
local maximum points, then for each positive integer $q$, there is 
a nonisotopic embedding
$\widetilde\Gamma_q$ of $\Gamma$ in which $C$ is replaced by a 
knot not isotopic to $C$, with $\nc(\widetilde\Gamma_q)$ as
close as desired to $\nc(p_e\circ\Gamma).$
\end{pro}
\pf
It follows from Corollary \ref{1dsuff} that the one-dimensional
graph $p_e \circ \Gamma$ may be replaced by an embedding
$\widehat\Gamma$ into a small neighborhood of the line $\Re e$ in 
$\Re^3$, with arbitrarily small change in its net total curvature.
Since $p_e \circ C$ has at least two local maximum points, there
is an interval 
of $\Re$ over which $p_e\circ C$ contains an interval which is the
image of 
four oriented intervals $J_1,J_2,J_3,J_4$ appearing in that cyclic
order around the oriented cycle $C$. Consider a plane presentation
of $\Gamma$ by orthogonal projection into a generic plane
containing the line $\Re e$.  Choose an integer $q\in\Z$, 
$|q|\geq 3.$ We modify $\widehat\Gamma$ by 
wrapping its interval $J_1$ $q$ times around $J_3$ 
and returning, passing over any other edges of $\Gamma$, including
$J_2$ and $J_4$, which it encounters along the way. The new graph
in $\Re^3$ is called $\widetilde{\Gamma_q}$. Then, 
if $C$ was the unknot, the cycle $\widetilde C_q$ which has
replaced it is a $(2,q)$-torus knot (see \cite{L}). In any case, 
$\widetilde C_q$ is not isotopic to $C$, and therefore
$\widetilde\Gamma_q$ is not isotopic to $\Gamma$.

As in the proof of Theorem \ref{incrdecr}, let 
$g_\delta:\Re^3 \to \Re^3$ be defined by cylindrical shrinking, so
that $g_1$ is the identity and $g_0=p_e$.
Then $p_e\circ \widetilde{\Gamma_q}=g_0(\widetilde{\Gamma_q})$, and
for $\delta>0$, $g_\delta(\widetilde{\Gamma_q})$ is isotopic to  
$\widetilde{\Gamma_q}$. But $\nc(g_\delta)\to \nc(g_0)$ as 
$\delta \to 0$.
\qed

\vspace{1em}

%
\begin{cor}\label{2cycle} 
If $e=e_0\in S^{n-1}$ minimizes $\nc(p_e\circ \Gamma)$, and there is a
cycle $C \subset \Gamma$ so that $p_{e_0}\circ C$ has 
two (or more) local maximum points, then there is a sequence of 
nonisotopic embeddings $\widetilde\Gamma_q$ of $\Gamma$ with
${\nc(\widetilde\Gamma_q)}$ less than, or 	
as close as desired, to $\nc(\Gamma)$,
in which $C$ is replaced by a $(2,q)$-torus knot.
\end{cor}

%
\begin{cor}\label{K_m2cycle}
If $\,\Gamma$ is an embedding of $K_m$ into
$\Re^3$, linear on each topological edge of $K_m$, $m\geq 4$,
then there is a sequence of nonisotopic embeddings
$\widetilde\Gamma_q$ of $\,\Gamma$ with
$\nc([\widetilde\Gamma_q])$ as
close as desired to $\nc([\Gamma])$, in which an unknotted cycle
$C$ of $\,\Gamma$ is replaced by a $(2,q)$-torus knot.
\end{cor}
\pf
According to Corollary \ref{2cycle}, we only need to construct an
isotopy of $K_m$ with the minimum value of $\nc$,
such that there is a cycle $C$ so that $p_e\circ C$ has two local
maximum points, where $\mu(e)$ is a minimum among $e \in S^2$. 

Choose $g:K_m\to\Re$ which is monotone on each edge of $K_m$, and 
has distinct values at vertices. Then according to Proposition
\ref{NTCK_m}, we have $\nc(g)=\nc(\{K_m\})$.  Number the vertices
$v_1,\dots,v_m$ so that $g(v_1)<g(v_2)<\dots<g(v_m)$.  Write
$E_{ji}$ for the edge  $E_{ij}$ with the reverse orientation, 
$i\neq j$.  Then the cycle $C$ formed in sequence from 
$E_{13},E_{32},E_{24}$ and $E_{41}$ has local maximum points at 
$v_3$ and $v_4$, and covers the interval
$\Big(g(v_2),g(v_3)\Big)\subset\Re$ four times. Since $C$ is
formed out of four straight edges, it is unknotted. The procedure
of Corollary \ref{2cycle} replaces $C$ with a $(2,q)$-torus knot,
with an arbitrarily small increase in NTC.
\qed \\

Note that 
Corollary \ref{2cycle} gives a set of conditions for those graph
types where a F\'{a}ry-Milnor type isotopy classification might 
hold. In particular, we consider one of the simpler
homeomorphism types of graphs, the 
{\bf theta graph}, $\theta=\theta_3=K_{3,2}$ (cf. description	
following Definition \ref{bridge}).  
The {\bf standard theta graph} is
the isotopy class in $\Re^3$ of a plane circle plus a diameter.
We have seen in Corollary \ref{theta_m} that the minimum of 
net total curvature for a theta graph is $3\pi$. On the other hand
note that in the range $3\pi\leq\nc(\Gamma) < 4\pi$, 
for $e$ in a set of positive measure of $S^2$, $p_e(\Gamma)$ 
cannot have two local maximum points. In Theorem \ref{thetathm}
below, we shall show that a theta graph $\Gamma$ with
$\nc(\Gamma)< 4\pi$ is isotopically standard.

We may observe that there are nonstandard theta graphs in
$\Re^3$.  For example, the union of two edges might form a knot.
Moreover, as S. Kinoshita has shown, there are 
$\theta$-graphs in $\Re^3$, not isotopic to a planar graph, such
that each of the three cycles formed by deleting one edge is
unknotted \cite{Ki}.

We begin with a well-known property of knots, whose proof we give
for the sake of completeness.

%
\begin{lem}\label{jordan}
Let $C \subset \Re^3$ be homeomorphic to $S^1$, and {\bf not} a
convex planar curve. Then there is a nonempty open set of planes 
$P\subset \Re^3$ which each meet $C$ in at least four points.
\end{lem}
\pf
For $e\in S^2$ and $t\in\Re$ write the plane
$P_t^e=\{x\in\Re^3:\langle e,x\rangle=t\}$.

If $C$ is not planar, then there exist four non-coplanar
points $p_1, p_2, p_3, p_4$, numbered in order around $C$.
Note that no three of the points can be collinear.  Let an
oriented plane $P_0$ be chosen to contain $p_1$ and $p_3$ and
rotated until both $p_2$ and $p_4$ are above $P_0$ strictly.
Write $e_1$ for the unit normal vector to $P_0$ on the side where
$p_2$ and $p_3$ lie, so that $P_0=P_{t_0=0}^{e_1}$. 	.
Then the set $P_t \cap C$ contains at least four points, for
$t_0=0<t<\delta_1$, with some $\delta_1>0$, since each plane
$P_t=P_t^{e_1}$ meets each of the four open arcs between the
points $p_1, p_2, p_3, p_4$.  This conclusion remains true, for
some $0<\delta < \delta_1$,  when the normal vector $e_1$ to $P_0$ is
replaced by any nearby $e \in S^2$, and $t$ is replaced by any
$0<t<\delta$.

If $C$ is planar but nonconvex, then there exists a plane 
$P_0=P_0^{e_1}$, transverse to the plane containing $C$, which
supports $C$ and touches $C$ at two distinct points,
but does not include the arc of $C$ between these two points.  
Consider disjoint open arcs of $C$ on either side of these
two points and including points not in $P_0$.  Then for 
$0 < t < \delta \ll 1$, the set $P_t \cap C$ contains at 
least four points, since the planes $P_t=P_t^{e_1}$ meet each 
of the four disjoint arcs.  Here once again $e_1$ may be 
replaced by any nearby unit vector $e$,  and the plane $P_t^e$
will meet $C$ in at least four points, for $t$ in a 
nonempty open interval $t_1<t<t_1+\delta$.
\qed\\

Using the notion of net total curvature, we may extend the
theorems of Fenchel \cite{Fen} as well as of F\'ary-Milnor
(\cite{Fa},\cite{M}), for curves homeomorpic to $S^1$, to graphs
homeomorphic to the theta graph. An analogous result is given by
Taniyama in \cite{T}, who showed that the minimum of $\tc$ for
polygonal 
$\theta$-graphs is $4\pi$, and that any $\theta$-graph $\Gamma$
with $\tc(\Gamma)<5\pi$ is isotopically standard, 

%
\begin{thm}\label{thetathm}
Suppose $f:\theta \to\Re^3$ is a continuous
embedding, $\Gamma=f(\theta)$.
Then $\nc(\Gamma) \geq 3\pi$.  If $\nc(\Gamma) < 4\pi$, then
$\Gamma$ is isotopic in $\Re^3$ to the planar theta graph.
Moreover, $\nc(\Gamma) = 3\pi$ iff the graph is a planar 
convex curve plus a straight chord.
\end{thm}

\pf
We consider first the case when $f:\theta\to\Re^3$ is piecewise
$C^2$. 

{\bf (1)} We have
shown the 
{\bf lower bound} $3\pi$ for $\nc(f)$, where $f:\theta\to \Re^n$ 
is any piecewise $C^2$ mapping,
since $\theta=\theta_3$ is 
one case of Corollary \ref{theta_m}, with $m=3$.

{\bf (2)} We show next that if there is a cycle $C$ in a graph $\Gamma$
(a subgraph homeomorphic to $S^1$)	
which satisfies the conclusion of Lemma \ref{jordan}, then 
$\mu(e) \geq 2$ for $e$ in a nonempty open set of $S^2$.  Namely,
for $t_0<t<t_0+\delta$, a family of planes $P_t^e$ meets $C$, and
therefore meets $\Gamma$, in at least four points. This is
equivalent to saying that the cardinality $\#(e,t)\geq 4$. This
implies, by Corollary \ref{fibcard}, that 
$\sum\{{\rm  nlm}(e,q): p_e(q) > t_0\}\geq 2$.  Thus,
since ${\rm nlm}^+(e,q)\geq {\rm nlm}(e,q)$, using Definition
\ref{defmu}, we have $\mu(e) \geq 2$.

Now consider the {\bf equality} case of a theta graph $\Gamma$
with $\nc(\Gamma) = 3\pi$.  As we have seen in the proof of
Proposition \ref{bipartite} with $m=3$ and $n=2$, the multiplicity 
$\mu(e) \geq \frac32=\frac{mn}{4}$ for a.a. $e \in S^2$, while the 
integral of
$\mu(e)$ over $S^2$ equals $ 2\,\nc(\Gamma) = 6\pi$ by Theorem
\ref{muthm}, implying $\mu(e) = 3/2$ a.e. on $S^2$. 
Thus, the conclusion of Lemma \ref{jordan} is impossible 
for any cycle $C$ in $\Gamma$. By Lemma \ref{jordan}, all 
cycles $C$ of $\Gamma$ must be planar and convex.

Now $\Gamma$ consists of three arcs $a_1$, $a_2$ and
$a_3$, 
with common endpoints $q^+$ and $q^-$. 
As we have just shown, the three Jordan curves 
$\Gamma_1:=a_2\cup a_3$,
$\Gamma_2:=a_3\cup a_1$ and $\Gamma_3:=a_1\cup a_2$ 	
are each planar and convex.  It follows that $\Gamma_1,\,\Gamma_2$
and $\Gamma_3$ lie in a common plane.  In terms of the topology of
this plane, one of the three arcs $a_1$, $a_2$ and $a_3$ lies in
the middle between the other two. But the middle arc, say $a_2$, 
must be a line segment, as it needs to be a shared piece of two 
curves $\Gamma_1$ and $\Gamma_3$ bounding disjoint convex open
sets in the plane.  The conclusion is that $\Gamma$ is a planar,
convex Jordan curve $\Gamma_2$, plus a straight chord $a_2$, 
whenever $\nc(\Gamma)= 3\pi.$

{\bf (3)} We next turn our attention to the 
{\bf upper bound} of $\nc$, to imply that a $\theta$-graph 
is isotopically standard:
we shall assume that $g:\theta \to \Re^3$ is an embedding in
general position with $\nc(g) < 4\pi$, 
and write $\Gamma=g(\theta)$.  	
By Theorem \ref{muthm}, since $S^2$ has area $4\pi$, the average
of $\mu(e)$ over $S^2$ is less than $2$, and it follows that there
exists a set of positive measure of $e_0\in S^2$ with 
$\mu(e_0) < 2$.  Since $\mu(e_0)$ is a half-integer, and since
$\mu(e) \geq 3/2$, as we have shown in 
part {\bf (1)} of this proof, we have $\mu(e_0) = 3/2$ exactly.

From Corollary \ref{muformula} applied to 
$p_{e_0}\circ g:\theta \to \Re$, we find 
$\mu_g(e_0)=\frac12(\Lambda+V)+\frac{k}{4}$, where $\Lambda$ is
the number of local maximum points, $V$ is the number of local
minimum points and $k=2$ is the number of vertices, both of
degree $3$. Thus, $\frac{3}{2}=\frac12(\Lambda+V)+\frac12$, so
that $\Lambda+V=2$. This implies that the local maximum/minimum
points are unique, and must be the unique global maximum/minimum
points $p_{\rm max}$ and $p_{\rm min}$ (which may be one of the
two vertices $q^\pm$).
Then $p_{e_0}\circ g$ is monotone along edges except at the 
points $p_{\rm max}$, $p_{\rm min}$ and $q^\pm$.

Introduce Euclidean coordinates $(x,y,z)$ for $\Re^3$ so that
$e_0$ is in the increasing $z$-direction. Write 
$t_{\rm max}=p_{e_0}\circ g(p_{\rm max})=
\langle e_0, p_{\rm max}\rangle$ and 
$t_{\rm min}=\langle e_0, p_{\rm min}\rangle$ for the maximum and
minimum values of $z$ along $g(\theta)$.  Write $t^\pm$ for the
value of $z$ at $g(q^\pm)$, where we may assume 
$t_{\rm min} \leq t^- < t^+ \leq t_{\rm max}$.

We construct a ``model" standard $\theta$-curve $\widehat\Gamma$ 
in the
$(x,z)$-plane, as follows. $\widehat\Gamma$ will consist of a
circle $C$ plus the straight chord of $C$,
joining $\widehat{q}^-$ to $\widehat{q}^+$ (points to be chosen). 
Choose $C$ so that the maximum and minimum values of $z$ on $C$
equal $t_{\rm max}$ and $t_{\rm min}$. 
Write $\widehat{p}_{\rm max}$ resp.
$\widehat{p}_{\rm min}$ for the maximum and minimum points of
$z$ along $C$. Choose $\widehat{q}^+$ as a point on
$C$ where $z=t^+$. There may be two nonequivalent choices for
$\widehat{q}^-$ as a point on $C$ where
$z=t^-$: we choose so that $\widehat{p}_{\rm max}$ and 
$\widehat{p}_{\rm min}$ are in the same or different topological
edge of $\widehat\Gamma$, where $p_{\rm max}$ and
$p_{\rm min}$ are in the same or different topological
edge, resp., of $\Gamma$. 
Note that there is a homeomorphism from
$g(\theta)$ to $\widehat\Gamma$ which preserves $z$. 

We now proceed to extend this homeomorphism to an isotopy. For
$t\in\Re$, write $P_t$ for the plane 
$\{z=t\}$.  As in the
proof of Proposition \ref{untangle}, there is a continuous
$1$-parameter family of homeomorphisms $\Phi_t:P_t\to P_t$ such
that $\Phi_t(\Gamma\cap P_t)=\widehat\Gamma\cap P_t$;
$\Phi_t$ is the identity outside a compact subset of $P_t$;
and $\Phi_t$ is isotopic to the identity of $P_t$, uniformly with
respect to $t$. Defining 
$\Phi:\Re^3\to\Re^3$ by $\Phi(x,y,z):=\Phi_z(x,y)$, we have an
isotopy of $\Gamma$ with the model graph $\widehat\Gamma$.

{\bf (4)} Finally, consider an embedding $g:\theta\to \Re^3$ which 
is only {\bf continuous}, and write $\Gamma=g(\theta)$. 

It follows from Theorem \ref{tame} that for any $\theta$-graph
$\Gamma$ of finite net total curvature, there is a
$\Gamma$-approximating polygonal $\theta$-graph $P$ isotopic to
$\Gamma$, with $\nc(P) \leq \nc(\Gamma)$ and as close as desired
to $\nc(\Gamma)$.

If a $\theta$-graph 
$\Gamma$ would have $\nc(\Gamma) < 3\pi$, then the
$\Gamma$-approximating
polygonal graph $P$ would also have $\nc(P) < 3\pi$, in
contradiction to what we have shown for piecewise $C^2$ theta
graphs in part {\bf (1)} above.
This shows that $\nc(\Gamma) \geq 3\pi$.

If equality $\nc(\Gamma) = 3\pi$ holds, then 
$\nc(P) \leq \nc(\Gamma) = 3\pi$,
so that by 
the equality case part {\bf (2)} above, $\nc(P)$ must equal
$3\pi$, and $P$ must be a convex planar curve plus a chord.  But
this holds for {\em all} $\Gamma$-approximating polygonal graphs
$P$, implying that 
$\Gamma$ itself must be a convex planar curve plus a chord.

Finally, If $\nc(\Gamma) < 4\pi$, then $\nc(P) < 4\pi$, implying
by part {\bf (3)}  above that $P$ is isotopic to the standard
$\theta$-graph.  But $\Gamma$ is isotopic to $P$, and hence is
isotopically standard.
\qed
\\

\bigskip

\normalsize

\begin{tabbing}
aaaaaaaaaaaaaaaaaaaaaaaaaaaaaasssssssssssssssss \=
bbbbbbbbbbbbbbbbbbbbbbbbbbbbbbb\kill

Robert Gulliver\> Sumio Yamada\\
School of Mathematics \>Mathematical Institute\\
University of Minnesota \> Tohoku University\\
Minneapolis MN 55414 \> Aoba, Sendai, Japan 980-8578\\
{\tt gulliver@math.umn.edu} \> {\tt yamada@math.tohoku.ac.jp}\\
{\tt www.math.umn.edu/\~{ }gulliver}\> {\tt
www.math.tohoku.ac.jp/\~{ }yamada}
\\
\end{tabbing}


\begin{thebibliography}{20}

\bibitem[AA]{AA} W. Allard and F. Almgren, {\it The structure of
stationary one dimensional varifolds with positive density},
Invent. Math {\bf 34} (1976), 83--97. 

\bibitem[AF]{AF} E. Artin and R. H. Fox, {\it Some wild cells and
spheres in three-dimensional space}, Annals of Math. {\bf 49}
(1948), 979-990. 

\bibitem[D1]{D1} J. Douglas, {\it  Solution of the problem of
Plateau}, Trans. Amer. Math. Soc. {\bf 33}(1931), 263--321.

\bibitem[EWW]{EWW} T. Ekholm, B. White, and D. Wienholtz, {\it
Embeddedness of minimal surfaces with total boundary curvature at
most $4\pi$}, Annals of Mathematics {\bf 155} (2002), 109--234.

\bibitem[Fa]{Fa}  I. F\'{a}ry,  {\it Sur la courbure totale
d'une courbe gauche faisant un noeud}, Bull. Soc. Math. France
{\bf 77} (1949), 128-138.

\bibitem[Fen]{Fen} W. Fenchel, {\it \"Uber Kr\"ummung und Windung
geschlossener Raumkurven}, Math. Ann. {\bf 101} (1929), 238--252.

\bibitem[G]{G} R. Gulliver, {\it Total Curvature of Graphs in
Space}, Pure and Applied Mathematics Quarterly  {\bf 3} (2007),
773--783.

\bibitem[GY1]{GY1} R. Gulliver and S. Yamada, {\it Area density
and regularity for soap film-like surfaces spanning graphs },
Math. Z. {\bf 253} (2006), 315--331.

\bibitem[GY2]{GY2} R. Gulliver and S. Yamada, {\it Total Curvature
and isotopy of graphs in $R^3$}, ArXiv:0806.0406.

\bibitem[Go]{Go} H. Goda, {\it Bridge index for theta curves in
the 3-sphere}, Topology Applic. {\bf 79} (1997), 177-196.

\bibitem[Ki]{Ki} S. Kinoshita, {\it On elementary ideals of
polyhedra in the $3$-sphere}, Pacific J. Math. {\bf 42} (1972),
89-98.

\bibitem[L]{L} Lickorish, W.B.R., {\it An Introduction to Knot
Theory.} Graduate Texts in Mathematics {\bf 175}. Springer, 1997.
Springer, 1971.

\bibitem[M]{M} J. Milnor, {\it On the total curvature of knots},
Annals of Math. {\bf 52} (1950), 248--257.

\bibitem[R]{R} T. Rad\'{o}, {\it On the Problem of Plateau}.
Springer, 1971.

\bibitem[vR]{vR} A. C. M. van Rooij, {\it The total curvature of
curves}, Duke Math. J.  {\bf 32}, 313--324 (1965).

\bibitem[S]{S} K. Stephenson, {Circle packing: a mathematical
tale} Notices AMS, {\bf 50} No.11, 1376--1388 (2003).

\bibitem[T]{T} K. Taniyama, {\it Total curvature of graphs in
Euclidean spaces}.  Differential Geom. Appl. {\bf 8} (1998),
135--155.

\end{thebibliography}
\end{document}